\documentclass[12pt,english]{article}
\usepackage[T1]{fontenc}
\usepackage[latin1]{inputenc}
\usepackage{amsmath}
\usepackage{amssymb}

\makeatletter

\usepackage{amsfonts}

\usepackage{amscd}

\usepackage{amsfonts}

\providecommand{\U}[1]{\protect\rule{.1in}{.1in}}
\newtheorem{theorem}{Theorem}[section]
\newtheorem{lemma}[theorem]{Lemma}
\newtheorem{proposition}[theorem]{Proposition}

\newtheorem{example}[theorem]{Example}
\newtheorem{definition}[theorem]{Definition}
\newtheorem{corollary}[theorem]{Corollary}
\newtheorem{remark}[theorem]{Remark}
\newenvironment{proof}{\bf Proof. \rm}{$\Box$}

\usepackage{babel}
\makeatother
\begin{document}

\title{Schur Class Operator Functions and Automorphisms of Hardy Algebras}

\author{Paul S. Muhly%
\thanks{Supported in part by grants from the National Science Foundation and
from the U.S.-Israel Binational Science Foundation.%
}\\
Department of Mathematics\\
University of Iowa\\
Iowa City, IA 52242\\
e-mail: pmuhly@math.uiowa.edu \and Baruch Solel%
\thanks{Supported in part by the U.S.-Israel Binational Science Foundation
and by the Fund for the Promotion of Research at the Technion.%
}\\
Department of Mathematics\\
Technion\\
32000 Haifa, Israel\\
e-mail: mabaruch@techunix.technion.ac.il}

\date{}

\maketitle
\begin{abstract}
Let $E$ be a $W^{\ast}$-correspondence over a von Neumann algebra
$M$ and let $H^{\infty}(E)$ be the associated Hardy algebra. If
$\sigma$ is a faithful normal representation of $M$ on a Hilbert
space $H$, then one may form the dual correspondence $E^{\sigma}$
and represent elements in $H^{\infty}(E)$ as $B(H)$-valued functions
on the unit ball $\mathbb{D}(E^{\sigma})^{\ast}$. The functions that
one obtains are called Schur class functions and may be characterized
in terms of certain Pick-like kernels. We study these functions and
relate them to system matrices and transfer functions from systems
theory. We use the information gained to describe the automorphism
group of $H^{\infty}(E)$ in terms of special Möbius transformations
on $\mathbb{D}(E^{\sigma})$. Particular attention is devoted to the
$H^{\infty}$-algebras that are associated to graphs.
\end{abstract}

\section{Introduction}

Let $M$ be a $W^{\ast}$-algebra and let $E$ be a $W^{\ast}$-correspondence
over $M$. In \cite{MS03} we built an operator algebra from this
data that we called the Hardy algebra of $E$ and which we denoted
$H^{\infty}(E)$. If $M=E=\mathbb{C}$ - the complex numbers, then
$H^{\infty}(E)$ is the classical Hardy algebra consisting of all
bounded analytic functions on the open unit disc, $\mathbb{D}$ (see
Example \ref{ex0.1} below.) If $M=\mathbb{C}$ again, but $E=\mathbb{C}^{n}$,
then $H^{\infty}(E)$ is the free semigroup algebra $\mathcal{L}_{n}$
studied by Davidson and Pitts \cite{DP98}, Popescu \cite{Po91} and
others (see Example \ref{ex0.2}.) One of the principal discoveries
made in \cite{MS03}, and the source of inspiration for the present
paper, is that attached to each faithful normal representation $\sigma$
of $M$ there is a \emph{dual} correspondence $E^{\sigma}$, which
is a $W^{\ast}$-correspondence over the commutant of $\sigma(M)$,
$\sigma(M)^{\prime}$, and the elements of $H^{\infty}(E)$ define
functions on the open unit ball of $E^{\sigma}$, $\mathbb{D}(E^{\sigma})$.
Further, the value distribution theory of these functions turns out
to be linked through our generalization of the Nevanlinna-Pick interpolation
theorem \cite[Theorem 5.3]{MS03} with the positivity properties of
certain Pick-like kernels of \emph{mappings between operator spaces}.

In the setting where $M=E=\mathbb{C}$ and $\sigma$ is the $1$-dimensional
representation of $\mathbb{C}$ on itself, then $E^{\sigma}$ is $\mathbb{C}$
again. The representation of $H^{\infty}(E)$ in terms of functions
on $\mathbb{D}{(E}^{\sigma})=\mathbb{D}$ is just the usual way we
think of $H^{\infty}(E)$. In this setting, our Nevanlinna-Pick theorem
is exactly the classical theorem. If, however, $\sigma$ is a representation
of $\mathbb{C}$ on a Hilbert space $H$, $\dim(H)>1$, then $E^{\sigma}$
may be identified with $B(H)$ and then $\mathbb{D}(E^{\sigma})$
becomes the space of strict contractions on $H$, i.e., all those
operators of norm strictly less than $1$. In this case, the value
of an $f\in H^{\infty}(E)$ at a $T\in\mathbb{D}{(E}^{\sigma})$ is
simply $f(T)$, defined through the usual holomorphic functional calculus.
Our Nevanlinna-Pick theorem gives a solution to problems such as this:
given $k$ operators $T_{1},T_{2},\ldots,T_{k}$ all of norm less
than $1$ and $k$ operators, $A_{1},A_{2},\ldots,A_{k}$, determine
the circumstances under which one can find a bounded analytic function
$f$ on the open unit disc of sup norm at most $1$ such that $f(T_{i})=A_{i}$,
$i=1,2,\ldots,k$ (See \cite[Theorem 6.1]{MS03}.) On the other hand,
when $M=\mathbb{C}$, $E=\mathbb{C}^{n}$, and $\sigma$ is one dimensional,
the space $E^{\sigma}$ is $\mathbb{C}^{n}$ and $\mathbb{D}(E^{\sigma})\mathbb{\ }$
is the unit ball $\mathbb{B}^{n}$. Elements in $H^{\infty}(E)=\mathcal{L}_{n}$
are realized as holomorphic functions on $\mathbb{B}^{n}$ that lie
in a multiplier space studied in detail by Arveson \cite{wA98}. More
accurately, the functional representation of $H^{\infty}(E)=\mathcal{L}_{n}$
in terms of these functions expresses this space as a \emph{quotient}
of $H^{\infty}(E)=\mathcal{L}_{n}$. The Nevanlinna-Pick theorem of
\cite{MS03} contains those of Davidson and Pitts \cite{DP98b}, Popescu
\cite{gP98}, and Arias and Popescu \cite{AP00}, which deal with
interpolation problems for these spaces of functions (possibly tensored
with the bounded operators on an auxiliary Hilbert space). It also
contains some of the results of Constaninescu and Johnson in \cite{CJ03}
which treats elements of $\mathcal{L}_{n}$ as functions on the ball
of strict row contractions with values in the operators on a Hilbert
space. (See their Theorem 3.4 in particular.) This situation arises
when one takes $M=\mathbb{C}$ and $E=\mathbb{C}^{n}$, but takes
$\sigma$ to be scalar multiplication on an auxiliary Hilbert space.

Our objective in the present note is basically two fold. First, we
wish to identify those functions on $\mathbb{D}{(E}^{\sigma})$ that
arise from evaluating elements of $H^{\infty}(E)$. For this purpose,
we introduce a family of functions on $\mathbb{D}{(E}^{\sigma})$
that we call Schur class operator functions (see Definition \ref{schur}).
Roughly speaking, these functions are defined so that a Pick-like
kernel that one may attach to each one is completely positive definite
in the sense of Barreto, Bhat, Liebscher and Skeide \cite{BBLS}.
In Theorem \ref{realization} we use their Theorem 3.2.3 to give a
Kolmogorov-type representation of the kernel, from which we derive
an analogue of a unitary system matrix $\left(\begin{array}{cc}
A & B\\
C & D\end{array}\right)$ whose transfer function \[
A+B(I-L_{\eta}^{\ast}D)^{-1}L_{\eta}^{\ast}C\]
 turns out to be the given Schur class operator function. We then
prove in Theorem \ref{evaluationtheorem} that each such transfer
function arises by evaluating an element in $H^{\infty}(E)$ at points
of $\mathbb{D}(E^{\sigma})$ and conversely, each function in $H^{\infty}(E)$
has a representation in terms of a transfer function. The meaning
of the notation will be made precise below, but we use it here to
highlight the connection between our analysis and realization theory
as it comes from mathematical systems theory. The point to keep in
mind is that functions on $\mathbb{D}(E^{\sigma})$ that come from
elements of $H^{\infty}(E)$ are not, \emph{a priori}, analytic in
any ordinary sense and it is not at all clear what analytic features
they have. Our Theorems \ref{schur} and \ref{evaluationtheorem}
together with \cite[Theorem 5.3]{MS03} show that the Schur class
operator functions are precisely the functions one obtains when evaluating
functions in $H^{\infty}(E)$ (of norm at most $1$) at points of
$\mathbb{D}(E^{\sigma})$. The fact that each such function may be
realized as a transfer function exhibits a surprising level of analyticity
that is not evident in the definition of $H^{\infty}(E)$.

Our second objective is to connect the usual holomorphic properties
of $\mathbb{D}(E^{\sigma})$ with the automorphisms of $H^{\infty}(E)$.
As a space, $\mathbb{D}(E^{\sigma})$ is the unit ball of a $J^{\ast}$-triple
system. Consequently, every holomorphic automorphism of $\mathbb{D}(E^{\sigma})$
is the composition of a Möbius transformation and a linear isometry
\cite{lH74}. Each of these implements an automorphism of the algebra
of all bounded, \emph{complex-valued analytic} functions on $\mathbb{D}(E^{\sigma})$,
but in our setting only certain of them implement automorphisms of
$H^{\infty}(E)$ - those for which the Möbius part is determined by
a {}``central\textquotedblright\ element of $E^{\sigma}$ (see
Theorem \ref{automcomp}). Our proof requires the fact that the evaluation
of functions in $H^{\infty}(E)$ (of norm at most $1$) at points
of $\mathbb{D}(E^{\sigma})$ are precisely the Schur class operator
functions on $\mathbb{D(}E^{\sigma})$. Indeed, the whole analysis
is an intricate {}``point - counterpoint\textquotedblright\ interplay
among elements of $H^{\infty}(E)$, Schur class functions, transfer
functions and {}``classical\textquotedblright\ function theory
on $\mathbb{D}(E^{\sigma})$. In the last section, we apply our general
analysis of the automorphisms of $H^{\infty}(E)$ to the special case
of $H^{\infty}$-algebras coming from directed graphs.

In concluding this introduction, we want to note that a preprint of
the present paper was posted on the arXiv on June 27, 2006. Recently,
inspired in part by our preprint, Ball, Biswas, Fang and ter Horst
\cite{BBFtHp07} were able to realize the Fock space that we describe
here in terms of the theory of completely positive definite kernels
advanced by Barreto, Bhat, Liebscher and Skeide \cite{BBLS} that
we also use (See Section 3 and, in particular, the proof of Theorem
\ref{realization}.) The analysis of Ball \emph{et al.} makes additional
ties between the theory of abstract Hardy algebras that we develop
here and classical function theory on the unit disc.

\section{Preliminaries}

We start by introducing the basic definitions and constructions. We
shall follow Lance \cite{L94} for the general theory of Hilbert $C^{\ast}$-modules
that we shall use. Let $A$ be a $C^{\ast}$-algebra and $E$ be a
right module over $A$ endowed with a bi-additive map $\langle\cdot,\cdot\rangle:E\times E\rightarrow A$
(referred to as an $A$-valued inner product) such that, for $\xi,\eta\in E$
and $a\in A$, $\langle\xi,\eta a\rangle=\langle\xi,\eta\rangle a$,
$\langle\xi,\eta\rangle^{\ast}=\langle\eta,\xi\rangle$, and $\langle\xi,\xi\rangle\geq0$,
with $\langle\xi,\xi\rangle=0$ only when $\xi=0$. Also, $E$ is
assumed to be complete in the norm $\Vert\xi\Vert:=\Vert\langle\xi,\xi\rangle\Vert^{1/2}$.
We write $\mathcal{L}(E)$ for the space of continuous, adjointable,
$A$-module maps on $E$. It is known to be a $C^{\ast}$-algebra.
If $M$ is a von Neumann algebra and if $E$ is a Hilbert $C^{\ast}$-module
over $M$, then $E$ is said to be \emph{self-dual} in case every
continuous $M$-module map from $E$ to $M$ is given by an inner
product with an element of $E$. Let $A$ and $B$ be $C^{\ast}$-algebras.
A ${C^{\ast}}$\emph{-correspondence} \emph{from} $A$\emph{\ to}
$B$ is a Hilbert $C^{\ast}$-module $E$ over $B$ endowed with a
structure of a left module over $A$ via a nondegenerate $\ast$-homomorphism
$\varphi:A\rightarrow\mathcal{L}(E)$.

When dealing with a specific $C^{\ast}$-correspondence, $E$, from
a $C^{\ast}$-algebra $A$ to a $C^{\ast}$-algebra $B$, it will
be convenient sometimes to suppress the $\varphi$ in formulas involving
the left action and simply write $a\xi$ or $a\cdot\xi$ for $\varphi(a)\xi$.
\ This should cause no confusion in context.

If $E$ is a $C^{\ast}$-correspondence from $A$ to $B$ and if $F$
is a correspondence from $B$ to $C$, then the balanced tensor product,
$E\otimes_{B}F$ is an $A,C$-bimodule that carries the inner product
defined by the formula \[
\langle\xi_{1}\otimes\eta_{1},\xi_{2}\otimes\eta_{2}\rangle_{E\otimes_{B}F}:=\langle\eta_{1},\varphi(\langle\xi_{1},\xi_{2}\rangle_{E})\eta_{2}\rangle_{F}\]
 The Hausdorff completion of this bimodule is again denoted by $E\otimes_{B}F$.

In this paper we deal mostly with correspondences over von Neumann
algebras that satisfy some natural additional properties as indicated
in the following definition. (For examples and more details see \cite{MS03}).

\begin{definition} Let $M$ and $N$ be von Neumann algebras and
let $E$ be a Hilbert $C^{\ast}$-module over $N$. Then $E$ is called
a \emph{Hilbert} $W^{\ast}$\emph{-module} over $N$ in case $E$
is self-dual. The module $E$ is called a $W^{\ast}$\emph{-correspondence
from} $M$\emph{\ to} $N$\emph{\ }in case $E$ is a self-dual $C^{\ast}$-correspondence
from $M$ to $N$ such that the $\ast$-homomorphism $\varphi:M\rightarrow\mathcal{L}(E)$,
giving the left module structure on $E$, is normal. If $M=N$ we
shall say that $E$ is a $W^{\ast}$-correspondence \emph{over $M$}.
\end{definition}

We note that if $E$ is a Hilbert $W^{\ast}$-module over a von Neumann
algebra, then $\mathcal{L}(E)$ is not only a $C^{\ast}$-algebra,
but is also a $W^{\ast}$-algebra. Thus it makes sense to talk about
normal homomorphisms into $\mathcal{L}(E)$.

\begin{definition} \label{isomorph}An \emph{isomorphism} of a $W^{\ast}$-correspondence
$E_{1}$ over $M_{1}$ and a $W^{\ast}$-correspondence $E_{2}$ over
$M_{2}$ is a pair $(\sigma,\Psi)$ where $\sigma:M_{1}\rightarrow M_{2}$
is an isomorphism of von Neumann algebras, $\Psi:E_{1}\rightarrow E_{2}$
is a vector space isomorphism preserving the $\sigma$-topology and
for $e,f\in E_{1}$ and $a,b\in M_{1}$, we have $\Psi(aeb)=\sigma(a)\Psi(e)\sigma(b)$
and $\langle\Psi(e),\Psi(f)\rangle=\sigma(\langle e,f\rangle)$. \end{definition}

When considering the tensor product $E\otimes_{M}F$ of two $W^{\ast}$-correspondences,
one needs to take the closure of the $C^{\ast}$-tensor product in
the $\sigma$-topology of \cite{BDH88} in order to get a $W^{\ast}$-correspondence.
However, we will not distinguish notationally between the $C^{\ast}$-tensor
product and the $W^{\ast}$-tensor product. Note also that given a
$W^{\ast}$-correspondence $E$ over $M$ and a Hilbert space $H$
equipped with a normal representation $\sigma$ of $M$, we can form
the Hilbert space $E\otimes_{\sigma}H$ by defining $\langle\xi_{1}\otimes h_{1},\xi_{2}\otimes h_{2}\rangle=\langle h_{1},\sigma(\langle\xi_{1},\xi_{2}\rangle)h_{2}\rangle$.
Thus, $H$ is viewed as a correspondence from $M$ to $\mathbb{C}$
via $\sigma$ and $E\otimes_{\sigma}H$ is just the tensor product
of $E$ and $H$ as $W^{\ast}$-correspondences.

Note also that, given an operator $X\in\mathcal{L}(E)$ and an operator
$S\in\sigma(M)^{\prime}$, the map $\xi\otimes h\mapsto X\xi\otimes Sh$
defines a bounded operator on $E\otimes_{\sigma}H$ denoted by $X\otimes S$.
The representation of $\mathcal{L}(E)$ that results when one lets
$S=I$, is called the representation of \emph{}$\mathcal{L}(E)$ \emph{induced
by} $\sigma$ and is often denoted by $\sigma^{E}$. The composition,
$\sigma^{E}\circ\varphi$ is a representation of $M$ which we shall
also say is induced by $\sigma$, but we shall usually denote it by
$\varphi(\cdot)\otimes I$.

Observe that if $E$ is a $W^{\ast}$-correspondence over a von Neumann
algebra $M$, then we may form the tensor powers $E^{\otimes n}$,
$n\geq0$, where $E^{\otimes0}$ is simply $M$ viewed as the identity
correspondence over $M$, and we may form the $W^{\ast}$-direct sum
of the tensor powers, $\mathcal{F}(E):=E^{\otimes0}\oplus E^{\otimes1}\oplus E^{\otimes2}\oplus\cdots$
to obtain a $W^{\ast}$-correspondence over $M$ called the (\emph{full})
\emph{Fock space} over $E$. The actions of $M$ on the left and right
of $\mathcal{F}(E)$ are the diagonal actions and, when it is convenient
to do so, we make explicit the left action by writing $\varphi_{\infty}$
for it. \ That is, for $a\in M$, $\varphi_{\infty}(a):=diag\{ a,\varphi(a),\varphi^{(2)}(a),\varphi^{(3)}(a),\cdots\}$,
where for all $n$, $\varphi^{(n)}(a)(\xi_{1}\otimes\xi_{2}\otimes\cdots\xi_{n})=(\varphi(a)\xi_{1})\otimes\xi_{2}\otimes\cdots\xi_{n}$,
$\xi_{1}\otimes\xi_{2}\otimes\cdots\xi_{n}\in E^{\otimes n}$. The
\emph{tensor algebra} over $E$, denoted $\mathcal{T}_{+}(E)$, is
defined to be the \emph{norm-closed} subalgebra of $\mathcal{L}(\mathcal{F}(E))$
generated by $\varphi_{\infty}(M)$ and the \emph{creation operators}
$T_{\xi}$, $\xi\in E$, defined by the formula $T_{\xi}\eta=\xi\otimes\eta$,
$\eta\in\mathcal{F}(E)$. We refer the reader to \cite{MS98} for
the basic facts about $\mathcal{T}_{+}(E)$.

\begin{definition} \label{Hinfty} (\cite{MS03}) Given a $W^{\ast}$-correspondence
$E$ over the von Neumann algebra $M$, the ultraweak closure of the
tensor algebra of $E$, $\mathcal{T}_{+}(E)$, in $\mathcal{L}(\mathcal{F}(E))$,
is called the \emph{Hardy Algebra of} $E$, and is denoted $H^{\infty}(E)$.
\end{definition}

\begin{example} \label{ex0.1}If $M=E=\mathbb{C}$, then $\mathcal{F}(E)$
can be identified with $\ell^{2}(\mathbb{Z}_{+})$ or, through the
Fourier transform, $H^{2}(\mathbb{T})$. The tensor algebra then is
isomorphic to the disc algebra $A(\mathbb{D})$ viewed as multiplication
operators on $H^{2}(\mathbb{T})$ and the Hardy algebra is realized
as the classical Hardy algebra $H^{\infty}(\mathbb{T})$. \end{example}

\begin{example} \label{ex0.2}If $M=\mathbb{C}$ and $E=\mathbb{C}^{n}$,
then $\mathcal{F}(E)$ can be identified with the space $l_{2}(\mathbb{F}_{n}^{+})$,
where $\mathbb{F}_{n}^{+}$ is the free semigroup on $n$ generators.
The tensor algebra then is what Popescu refers to as the {}``non
commutative disc algebra\textquotedblright\ $\mathcal{A}_{n}$ and
the Hardy algebra is its $w^{\ast}$-closure. It was studied by Popescu
\cite{Po91} and by Davidson and Pitts who denoted it by $\mathcal{L}_{n}$
\cite{DP98}. \end{example}

We need to review some basic facts about the representation theory
of $H^{\infty}(E)$ and of $\mathcal{T}_{+}(E)$. See \cite{MS98,MS03}
for more details.

\begin{definition} \label{Definition1.12}Let $E$ be a $W^{\ast}$-correspondence
over a von Neumann algebra $M$. Then:

\begin{enumerate}
\item A \emph{completely contractive covariant representation} of $E$ on
a Hilbert space $H$ is a pair $(T,\sigma)$, where

\begin{enumerate}
\item $\sigma$ is a normal $\ast$-representation of $M$ in $B(H)$.
\item $T$ is a linear, completely contractive map from $E$ to $B(H)$
that is continuous in the $\sigma$-topology of \cite{BDH88} on $E$
and the ultraweak topology on $B(H).$
\item $T$ is a bimodule map in the sense that $T(S\xi R)=\sigma(S)T(\xi)\sigma(R)$,
$\xi\in E$, and $S,R\in M$. 
\end{enumerate}
\item A completely contractive covariant representation $(T,\sigma)$ of
$E$ in $B(H)$ is called \emph{isometric} in case \begin{equation}
T(\xi)^{\ast}T(\eta)=\sigma(\langle\xi,\eta\rangle)\label{isometric}\end{equation}
 for all $\xi,\eta\in E$. 
\end{enumerate}
\end{definition}

It should be noted that the operator space structure on $E$ to which
Definition \ref{Definition1.12} refers is that which $E$ inherits
when viewed as a subspace of its linking algebra.

As we showed in \cite[Lemmas 3.4--3.6]{MS98} and in \cite{MS03},
if a completely contractive covariant representation, $(T,\sigma)$,
of $E$ in $B(H)$ is given, then it determines a contraction $\tilde{T}:E\otimes_{\sigma}H\rightarrow H$
defined by the formula $\tilde{T}(\eta\otimes h):=T(\eta)h$, $\eta\otimes h\in E\otimes_{\sigma}H$.
The operator $\tilde{T}$ intertwines the representation $\sigma$
on $H$ and the induced representation $\sigma^{E}\circ\varphi=\varphi(\cdot)\otimes I_{H}$
on $E\otimes_{\sigma}H$; i.e. \begin{equation}
\tilde{T}(\varphi(\cdot)\otimes I)=\sigma(\cdot)\tilde{T}.\label{covariance}\end{equation}
 In fact we have the following lemma from \cite[Lemma 2.16]{MS03}.

\begin{lemma} \label{CovRep}The map $(T,\sigma)\rightarrow\tilde{T}$
is a bijection between all completely contractive covariant representations
$(T,\sigma)$ of $E$ on the Hilbert space $H$ and contractive operators
$\tilde{T}:E\otimes_{\sigma}H\rightarrow H$ that satisfy equation
(\ref{covariance}). Given such a $\tilde{T}$ satisfying this equation,
$T$, defined by the formula $T(\xi)h:=\tilde{T}(\xi\otimes h)$,
together with $\sigma$ is a completely contractive covariant representation
of $E$ on $H$. Further, $(T,\sigma)$ is isometric if and only if
$\tilde{T}$ is an isometry. \end{lemma}

The importance of the completely contractive covariant representations
of $E$ (or, equivalently, the intertwining contractions $\tilde{T}$
as above) is that they yield all completely contractive representations
of the tensor algebra. More precisely, we have the following.

\begin{theorem} \label{Theorem310MS98}Let $E$ be a $W^{\ast}$-correspondence
over a von Neumann algebra $M$. To every completely contractive covariant
representation, $(T,\sigma)$, of $E$ there is a unique completely
contractive representation $\rho$ of the tensor algebra $\mathcal{T}_{+}(E)$
that satisfies \[
\rho(T_{\xi})=T(\xi)\;\;\;\xi\in E\]
 and\[
\rho(\varphi_{\infty}(a))=\sigma(a)\;\;\; a\in M.\]
 The map $(T,\sigma)\mapsto\rho$ is a bijection between the set of
all completely contractive covariant representations of $E$ and all
completely contractive (algebra) representations of $\mathcal{T}_{+}(E)$
whose restrictions to $\varphi_{\infty}(M)$ are continuous with respect
to the ultraweak topology on $\mathcal{L}(\mathcal{F}(E))$. \end{theorem}

\begin{definition} \label{integratedform}If $(T,\sigma)$ is a completely
contractive covariant representation of a $W^{\ast}$-correspondence
$E$ over a von Neumann algebra $M$, we call the representation $\rho$
of $\mathcal{T}_{+}(E)$ described in Theorem \ref{Theorem310MS98}
the \emph{integrated form} of $(T,\sigma)$ and write $\rho=\sigma\times T$.
\end{definition}

\begin{remark} \label{keyproblem}One of the principal difficulties
one faces in dealing with $\mathcal{T}_{+}(E)$ and $H^{\infty}(E)$
is to decide when the integrated form, $\sigma\times T$, of a completely
contractive covariant representation $(T,\sigma)$ extends from $\mathcal{T}_{+}(E)$
to $H^{\infty}(E)$. This problem arises already in the simplest situation,
vis. when $M=\mathbb{C}=E$. In this setting, $T$ is given by a single
contraction operator on a Hilbert space, $\mathcal{T}_{+}(E)$ {}``is\textquotedblright\ the
disc algebra and $H^{\infty}(E)$ {}``is\textquotedblright\ the
space of bounded analytic functions on the disc. The representation
$\sigma\times T$ extends from the disc algebra to $H^{\infty}(E)$
precisely when there is no singular part to the spectral measure of
the minimal unitary dilation of $T$. We are not aware of a comparable
result in our general context but we have some sufficient conditions.
One of them is given in the following lemma. It is not a necessary
condition in general. \end{remark}

\begin{lemma} \label{contraction} \cite[Corollary 2.14]{MS03} If
$\Vert\tilde{T}\Vert<1$ then $\sigma\times T$ extends to a ultraweakly
continuous representation of $H^{\infty}(E)$. \end{lemma}

In \cite{MS03} we introduced and studied the concepts of duality
and of point evaluation (for elements of $H^{\infty}(E)$). These
play a central role in our analysis here.

\begin{definition} \label{dual}Let $E$ be a $W^{\ast}$-correspondence
over a von Neumann algebra $M$ and let $\sigma:M\rightarrow B(H)$
be a faithful normal representation of $M$ on a Hilbert space $H$.
Then the $\sigma$\emph{-dual} of $E$, denoted $E^{\sigma}$, is
defined to be \[
\{\eta\in B(H,E\otimes_{\sigma}H)\mid\eta\sigma(a)=(\varphi(a)\otimes I)\eta,\; a\in M\}.\]

\end{definition}

An important feature of the dual $E^{\sigma}$ is that it is a $W^{\ast}$-correspondence,
\emph{but over the commutant} of $\sigma(M)$, $\sigma(M)^{\prime}$.

\begin{proposition} \label{corres} With respect to the action of
$\sigma(M)^{\prime}$ and the $\sigma(M)^{\prime}$-valued inner product
defined as follows, $E^{\sigma}$ becomes a $W^{\ast}$-correspondence
over $\sigma(M)^{\prime}$: For $Y$ and $X$ in $\sigma(M)^{\prime}$,
and $\eta\in E^{\sigma}$, $X\cdot\eta\cdot Y:=(I\otimes X)\eta Y$,
and for $\eta_{1},\eta_{2}\in E^{\sigma}$, $\langle\eta_{1},\eta_{2}\rangle_{\sigma(M)^{\prime}}:=\eta_{1}^{\ast}\eta_{2}$.
\end{proposition}

In the following remark we explain what we mean by {}``evaluating
an element of $H^{\infty}(E)$ at a point in the open unit ball of
the dual''.

\begin{remark} \label{evaluation}The importance of this dual space,
$E^{\sigma}$, is that it is closely related to the representations
of $E$. In fact, the operators in $E^{\sigma}$ whose norm does not
exceed $1$ are precisely the adjoints of the operators of the form
$\tilde{T}$ for a covariant pair $(T,\sigma)$. In particular, every
$\eta$ in the \emph{open} unit ball of $E^{\sigma}$ (written $\mathbb{D}(E^{\sigma})$)
gives rise to a covariant pair $(T,\sigma)$ (with $\eta=\tilde{T}^{\ast}$)
such that $\sigma\times T$ extends to a representation of $H^{\infty}(E)$.

Given $X\in H^{\infty}(E)$ we can apply the representation associated
to $\eta$ to it. The resulting operator in $B(H)$ will be denoted
by $\widehat{X}(\eta^{\ast})$. Thus \[
\widehat{X}(\eta^{\ast})=(\sigma\times\eta^{*})(X).\]
 In this way, we view every element in the Hardy algebra as a $B(H)$-valued
function \[
\widehat{X}:\mathbb{D}(E^{\sigma})^{*}\rightarrow B(H)\]
 on the open unit ball of $(E^{\sigma})^{*}$. One of our primary
objectives is to understand the range of the transform $X\rightarrow\widehat{X}\,$,
$X\in H^{\infty}(E)$. \end{remark}

\begin{example} \label{ex0.12}Suppose $M=E=\mathbb{C}$ and $\sigma$
the representation of $\mathbb{C}$ on some Hilbert space $H$. Then
it is easy to check that $E^{\sigma}$ is isomorphic to $B(H)$. Fix
an $X\in H^{\infty}(E)$. As we mentioned above, this Hardy algebra
is the classical $H^{\infty}(\mathbb{T})$ and we can identify $X$
with a function $f\in H^{\infty}(\mathbb{T})$. Given $S\in\mathbb{D}(E^{\sigma})=B(H)$,
it is not hard to check that $\widehat{X}(S^{\ast})$, as defined
above, is the operator $f(S^{\ast})$ defined through the usual holomorphic
functional calculus. \end{example}

\begin{example} \label{freesgp}In \cite{DP98} Davidson and Pitts
associate to every element of the free semigroup algebra $\mathcal{L}_{n}$
(see Example~\ref{ex0.2}) a function on the open unit ball of $\mathbb{C}^{n}$.
This is a special case of our analysis when $M=\mathbb{C}$, $E=\mathbb{C}^{n}$
and $\sigma$ is a one dimensional representation of $\mathbb{C}$.
In this case $\sigma(M)^{\prime}=\mathbb{C}$ and $E^{\sigma}=\mathbb{C}^{n}$.
Note, however, that our definition allows us to take $\sigma$ to
be the representation of $\mathbb{C}$ on an arbitrary Hilbert space
$H$. If we do so, then $E^{\sigma}$ is isomorphic to $B(H)^{(n)}$,
the nth column space over $B(H)$, and elements of $\mathcal{L}_{n}$
define functions on the open unit ball of this space viewed as a correspondence
over $B(H)$ with values in $B(H)$. This is the perspective adopted
by Constantinescu and Johnson in \cite{CJ03}. In the analysis of
\cite{DP98} it is possible that a non zero element of $\mathcal{L}_{n}$
will give rise to the zero function. We shall show in Lemma~\ref{1to1}
that, by choosing an appropriate $H$ we can insure that this does
not happen. \end{example}

\begin{example} \label{FreeHolFcns}Part of the recent work of Popescu
in \cite{gPp06} may be cast in our framework. We will follow his
notation. Fix a Hilbert space $K$, and let $E$ be the column space
$B(K)^{n}$. Take, also, a Hilbert space $H$ and let $\sigma:B(K)\rightarrow B(K\otimes H)$
be the representation which sends $a\in B(K)$ to $a\otimes I_{H}$.
Then, since the commutant of $\sigma(B(K))$ is naturally isomorphic
to $B(H)$, it is easy to see that $E^{\sigma}$ is the column space
over $B(H)$, $B(H)^{n}$. It follows that $\mathbb{D}(E^{\sigma})$
is the open unit ball in $B(H)^{n}$. A \emph{free formal power series}
with coefficients from $B(K)$ is a formal series $F=\sum_{\alpha\in\mathbb{F}_{n}^{+}}A_{\alpha}\otimes Z_{\alpha}$
where $\mathbb{F}_{n}^{+}$ is the free semigroup on $n$ generators,
the $A_{\alpha}$ are elements of $B(K)$ and where $Z_{\alpha}$
is the monomial in noncommuting indeterminates $Z_{1}$, $Z_{2}$,
\ldots{}, $Z_{n}$ determined by $\alpha$. If $F$ has radius of
convergence equal to $1$, then one may evaluate $F$ at points of
$\mathbb{D}(E^{\sigma})^{\ast}$ to get a function on $\mathbb{D}(E^{\sigma})^{\ast}$
with values in $B(K\otimes H)$, \emph{vis.,} $F((S_{1},S_{2},\cdots S_{n}))=\sum_{\alpha\in\mathbb{F}_{n}^{+}}A_{\alpha}\otimes S_{\alpha}$.
See \cite[Theorem 1.1]{gPp06}. In fact, under additional restrictions
on the coefficients $A_{\alpha}$, $F$ may be viewed as a function
$X$ in $H^{\infty}(B(K)^{n})$ in such a way that $F((S_{1},S_{2},\cdots S_{n}))=\widehat{X}(S_{1},S_{2},\cdots S_{n})$
in the sense defined in \cite[p. 384]{MS03} and discussed above in
Remark \ref{evaluation}. The space that Popescu denotes by $H^{\infty}(B(\mathcal{X)}_{1}^{n})$
arises when $K=\mathbb{C}$, and is naturally isometrically isomorphic
to $\mathcal{L}_{n}$ \cite[Theorem 3.1]{gPp06}. We noted in the
preceding example that $\mathcal{L}_{n}$ is $H^{\infty}(\mathbb{C}^{n})$.
The point of \cite{gPp06}, at least in part, is to study $H^{\infty}(B(\mathcal{X)}_{1}^{n})\simeq\mathcal{L}_{n}=H^{\infty}(\mathbb{C}^{n})$
through all the representations $\sigma$ of $\mathbb{C}$ on Hilbert
spaces $H$, that is, through evaluating functions in $H^{\infty}(B(\mathcal{X)}_{1}^{n})$
at points the unit ball of $B(H)^{n}$ for all possible $H$'s. The
space $B(K)^{n}$ is Morita equivalent to $\mathbb{C}^{n}$ in the
sense of \cite{MS00}, at least when $\dim(K)<\infty$, and, in that
case the tensor algebras $\mathcal{T}_{+}(B(K)^{n})$ and $\mathcal{T}_{+}(\mathbb{C}^{n})$
are Morita equivalent in the sense described by \cite{BMP00}. The
tensor algebra $\mathcal{T}_{+}(\mathbb{C}^{n})$, in turn, is naturally
isometrically isomorphic to Popescu's noncommutative disc algebra
$\mathcal{A}_{n}$ \cite{gP96}. The analysis in \cite{BMP00} suggests
a sense in which $\mathbb{C}^{n}$ and $B(K)^{n}$ are Morita equivalent
even when $\dim(K)=\infty$, and that together with \cite{MS00} suggests
that $H^{\infty}(B(K)^{n})$ should be Morita equivalent to $H^{\infty}(B(\mathcal{X)}_{1}^{n})\simeq H^{\infty}(\mathbb{C}^{n})$.
This would suggest an even closer connection between Popescu's free
power series, and all that goes with them, and the perspective we
have taken in this paper, which, as we shall see, involves generalized
Schur functions and transfer functions. The connection seems like
a promising avenue to explore. \end{example}

In \cite{MS03} we exploited the perspective of viewing elements of
the Hardy algebra as $B(H)$-valued functions on the open unit ball
of the dual correspondence to prove a Nevanlinna-Pick type interpolation
theorem. In order to state it we introduce some notation: For operators
$B_{1}$ and $B_{2}$ in $B(H)$, we write $Ad(B_{1},B_{2})$ for
the map from $B(H)$ to itself that sends $S$ to $B_{1}SB_{2}^{\ast}$.
Also, given elements $\eta_{1},\eta_{2}$ in $\mathbb{D}(E^{\sigma})$,
we let $\theta_{\eta_{1},\eta_{2}}$ denote the map, from $\sigma(M)^{\prime}$
to itself that sends $a$ to $\langle\eta_{1},a\eta_{2}\rangle$.
That is, $\theta_{\eta_{1},\eta_{2}}(a):=\langle\eta_{1},a\eta_{2}\rangle=\eta_{1}^{\ast}a\eta_{2}$,
$a\in\sigma(M)^{\prime}$.

\begin{theorem} \label{NP}(\cite[Theorem 5.3]{MS03}) Let $E$ be
a $W^{\ast}$-correspondence over a von Neumann algebra $M$ and let
$\sigma:M\rightarrow B(H)$ be a faithful normal representation of
$M$ on a Hilbert space $H$. Fix $k$ points $\eta_{1},\ldots\eta_{k}\ $in
the disk $\mathbb{D}(E^{\sigma})$ and choose $2k$ operators $B_{1},\ldots B_{k},C_{1},\ldots C_{k}$
in $B(H)$. Then there exists an $X$ in $H^{\infty}(E)$ such that
$\Vert X\Vert\leq1$ and\[
B_{i}\widehat{X}(\eta_{i}^{\ast})=C_{i}\]
 for $i=1,2,\ldots,k,$ if and only if the map from $M_{k}(\sigma(M)^{\prime})$
into $M_{k}(B(H))$ defined by the $k\times k$ matrix \begin{equation}
\left((Ad(B_{i},B_{j})-Ad(C_{i},C_{j}))\circ(id-\theta_{\eta_{i},\eta_{j}})^{-1}\right)\label{PickMatrix}\end{equation}
 is completely positive. \end{theorem}

That is, the map $T$, say, given by the matrix (\ref{PickMatrix})
is computed by the formula\[
T((a_{ij}))=(b_{ij}),\]
 where \[
b_{ij}=B_{i}((id-\theta_{\eta_{i},\eta_{j}})^{-1}(a_{ij})B_{j}^{\ast}-C_{i}((id-\theta_{\eta_{i},\eta_{j}})^{-1}(a_{ij})C_{j}^{\ast}\]
 and \[
(id-\theta_{\eta_{i},\eta_{j}})^{-1}(a_{ij})=a_{ij}+\theta_{\eta_{i},\eta_{j}}(a_{ij})+\theta_{\eta_{i},\eta_{j}}(\theta_{\eta_{i},\eta_{j}}(a_{ij}))+\cdots\]

We close this section with two technical lemmas that will be needed
in our analysis. Let $M$ and $N$ be $W^{\ast}$-algebras and let
$E$ be a $W^{\ast}$-correspondence from $M$ to $N$. Given a $\sigma$-closed
subcorrespondence $E_{0}$ of $E$ we know that the orthogonal projection
$P$ of $E$ onto $E_{0}$ is a right module map. (See \cite[Consequences 1.8 (ii)]{BDH88}).
In the following lemma we show that $P$ also preserves the left action.

\begin{lemma} \label{bimodule}Let $E$ be a $W^{\ast}$-correspondence
from the von Neumann algebra $M$ to the von Neumann algebra $N$,
and let $E_{0}$ be a sub $W^{\ast}$-correspondence $E_{0}$ of $E$
that is closed in the $\sigma$-topology of \cite[Consequences 1.8 (ii)]{BDH88}.
If $P$ is the orthogonal projection from $E$ onto $E_{0}$, then
$P$ is a bimodule map; i.e., $P(a\xi b)=aP(\xi)b$ for all $a\in M$
and $b\in N$. \end{lemma}

\begin{proof} It suffices to check that $P(e\xi)=eP(\xi)$ for all
$\xi\in E$ and projections $e\in M$. For $\xi,\eta\in E$ and a
projection $e\in M$, we have \[
\Vert e\xi+f\eta\Vert^{2}=\Vert\langle e\xi,e\xi\rangle+\langle f\eta,f\eta\rangle\Vert\leq\Vert\langle e\xi,e\xi\rangle\Vert+\Vert\langle f\eta,f\eta\rangle\Vert=\Vert e\xi\Vert^{2}+\Vert f\eta\Vert^{2},\]
 where $f=1-e$. So, for every $\lambda\in\mathbb{R}$ we have \[
(\lambda+1)^{2}\Vert fP(e\xi)\Vert^{2}=\Vert fP(e\xi+\lambda fP(e\xi))\Vert^{2}\leq\Vert e\xi+\lambda fP(e\xi)\Vert^{2}\]
\[
\leq\Vert e\xi\Vert^{2}+\lambda^{2}\Vert fP(e\xi)\Vert^{2}.\]
 Hence, for every $\lambda\in\mathbb{R}$, \[
(2\lambda+1)\Vert fP(e\xi)\Vert^{2}\leq\Vert e\xi\Vert^{2}\]
 and, thus, \[
(I-e)P(e\xi)=fP(e\xi)=0.\]
 Replacing $e$ by $f=I-e$ we get $eP((I-e)\xi)=0$ and, therefore,
\[
P(e\xi)=eP(e\xi)=eP(\xi).\]
 Since $M$ is spanned by its projections, we are done. \end{proof}

\begin{lemma} \label{faithful}Let $E$ be a $W^{\ast}$-correspondence
over $M$, let $\sigma$ be a faithful normal representation of $M$
on the Hilbert space $\mathcal{E}$, and let $E^{\sigma}$ be the
$\sigma$-dual correspondence over $N:=\sigma(M)^{\prime}$. Then

\begin{itemize}
\item [(i)] The left action of $N$ on $E^{\sigma}$ is faithful if and
only if $E$ is full (i.e. if and only if the ultraweakly closed ideal
generated by the inner products $\langle\xi_{1},\xi_{2}\rangle$,
$\xi_{1},\xi_{2}\in E$, is all of $M$).
\item [(ii)] The left action of $M$ on $E$ is faithful if and only if
$E^{\sigma}$ is full. 
\end{itemize}
\end{lemma}

\begin{proof} We shall prove (i). Part (ii) then follows by duality
(using \cite[Theorem 3.6]{MS03}). Given $S\in N$, $S\eta=0$ for
every $\eta\in E^{\sigma}$ if and only if for all $\eta\in E^{\sigma}$
and $g\in\mathcal{E}$, $(I\otimes S)\eta(g)=0$. Since the closed
subspace spanned by the ranges of all $\eta\in E^{\sigma}$ is all
of $E\otimes_{M}\mathcal{E}$ (\cite{MS03}), this is equivalent to
the equation $\xi\otimes Sg=0$ holding for all $g\in\mathcal{E}\;$and
$\xi\in E$. Since $\langle\xi\otimes Sg,\xi\otimes Sg\rangle=\langle g,S^{\ast}\langle\xi,\xi\rangle Sg\rangle$,
we find that $SE^{\sigma}=0$ if and only if $\sigma(\langle E,E\rangle)S=0$,
where $\langle E,E\rangle$ is the ultraweakly closed ideal generated
by all inner products. If this ideal is all of $M$ we find that the
equation $SE^{\sigma}=0$ implies that $S=0$. In the other direction,
if this is not the case, then this ideal is of the form $(I-q)M$
for some central nonzero projection $q$ and then $S=\sigma(q)$ is
different from $0$ but vanishes on $E^{\sigma}$. \end{proof}

\section{Schur class operator functions and realization}

Throughout this section, $E$ will be a fixed $W^{\ast}$-correspondence
over the von Neumann algebra $M$ and $\sigma$ will be a faithful
representation of $M$ on a Hilbert space $\mathcal{E}$. We then
form the $\sigma$-dual of $E$, $E^{\sigma}$, which is a correspondence
over $N:=\sigma(M)^{\prime}$, and we write $\mathbb{D}(E^{\sigma})$
for its open unit ball. Further, we write $\mathbb{D}(E^{\sigma})^{\ast}$
for $\{\eta^{\ast}\mid\eta\in\mathbb{D}(E^{\sigma})\}$. 

The following definition is clearly motivated by the condition appearing
in Theorem~\ref{NP} and Schur's theorem from classical function
theory.

\begin{definition} \label{schur}Let $\Omega$ be a subset of $\mathbb{D}(E^{\sigma})$
and let $\Omega^{\ast}=\{\omega^{\ast}\mid\omega\in\Omega\}$. A function
$Z:\Omega^{\ast}\rightarrow B(\mathcal{E})$ will be called a \emph{Schur
class operator function} (with values in $B(\mathcal{E})$) if, for
every $k$ and every choice of elements $\eta_{1},\eta_{2},\ldots,\eta_{k}$
in $\Omega$, the map from $M_{k}(N)$ to $M_{k}(B(\mathcal{E}))$
defined by the $k\times k$ matrix of maps, \[
((id-Ad(Z(\eta_{i}^{\ast}),Z(\eta_{j}^{\ast})))\circ(id-\theta_{\eta_{i},\eta_{j}})^{-1}),\]
 is completely positive. \end{definition}

Note that, when $M=E=B(\mathcal{E})$ and $\sigma$ is the identity
representation of $B(\mathcal{E})$ on $\mathcal{E}$, $\sigma(M)^{\prime}$
is $\mathbb{C}I_{\mathcal{E}}$, $E^{\sigma}$ is isomorphic to $\mathbb{C}$
and $\mathbb{D}(E^{\sigma})^{\ast}$ can be identified with the open
unit disc $\mathbb{D}$ of $\mathbb{C}$. In this case our definition
recovers the classical Schur class functions. More precisely, these
functions are usually defined as analytic functions $Z$ from an open
subset $\Omega$ of $\mathbb{D}$ into the closed unit ball of $B(\mathcal{E})$
but it is known that such functions are precisely those for which
the Pick kernel $k_{Z}(z,w)=(I-Z(z)Z(w)^{\ast})(1-z\bar{w})^{-1}$
is positive semi-definite on $\Omega$. The argument of \cite[Remark 5.4]{MS03}
shows that the positivity of this kernel is equivalent, in our case,
to the condition of Definition~\ref{schur}. This condition, in turn,
is the same as asserting that the kernel\begin{equation}
k_{Z}(\zeta^{\ast},\omega^{\ast}):=(id-Ad(Z(\zeta^{\ast}),Z(\omega^{\ast}))\circ(id-\theta_{\zeta,\omega})^{-1}\label{cpkernel}\end{equation}
 is a completely positive definite kernel on $\Omega^{\ast}$ in the
sense of Definition 3.2.2 of \cite{BBLS}.

For the sake of completeness, we record the fact that every element
of $H^{\infty}(E)$ of norm at most one gives rise to a Schur class
operator function.

\begin{theorem} \label{SchurHinfty}Let $E$ be a $W^{\ast}$-correspondence
over a von Neumann algebra $M$ and let $\sigma$ be a faithful normal
representation of $M$ in $B(H)$ for some Hilbert space $H$. If
$X$ is an element of $H^{\infty}(E)$ of norm at most one, then the
function $\eta^{\ast}\rightarrow\widehat{X}(\eta^{\ast})$ defined
in Remark \ref{evaluation} is a Schur class operator function on
$\mathbb{D}((E^{\sigma}))^{\ast}$ with values in $B(H)$. \end{theorem}

\begin{proof} One simply takes $B_{i}=I$ for all $i$ and $C_{i}=\widehat{X}(\eta_{i}^{\ast})$
in Theorem \ref{NP}. \end{proof}

\begin{theorem} \label{realization}Let $E$ be a $W^{\ast}$-correspondence
over a von Neumann algebra $M$. Suppose also that $\sigma$ a faithful
normal representation of $M$ on a Hilbert space $\mathcal{E}$ and
that $q_{1}$ and $q_{2}$ are projections in $\sigma(M)$. Finally,
suppose that $\Omega$ is a subset of $\mathbb{D}((E^{\sigma}))$
and that $Z$ is a Schur class operator function on $\Omega^{\ast}$
with values in $q_{2}B(\mathcal{E})q_{1}$. Then there is a Hilbert
space $H$, a normal representation $\tau$ of $N:=\sigma(M)^{\prime}$
on $H$ and operators $A,B,C$ and $D$ fulfilling the following conditions:

\begin{itemize}
\item [(i)] The operator $A$ lies in $q_{2}\sigma(M)q_{1}$.
\item [(ii)] The operators $C$, $B$, and $D$, are in the spaces $B(\mathcal{E}_{1},E^{\sigma}\otimes_{\tau}H)$,
$B(H,\mathcal{E}_{2})$, and $B(H,E^{\sigma}\otimes_{\tau}H)$, respectively,
and each intertwines the representations of $N=\sigma(M)^{\prime}$
on the relevant spaces (i.e. , for every $S\in N$, $CS=(S\otimes I_{H})C$,
$B\tau(S)=SB$ and $D\tau(S)=(S\otimes I_{H})D$).
\item [(iii)] The operator matrix \begin{equation}
V=\left(\begin{array}{cc}
A & B\\
C & D\end{array}\right),\label{transfermatrix}\end{equation}
 viewed as an operator from $\mathcal{E}_{1}\oplus H$ to $\mathcal{E}_{2}\oplus(E^{\sigma}\otimes_{\tau}H)$,
is a coisometry, which is unitary if $E$ is full.
\item [(iv)] For every $\eta^{\ast}$ in $\Omega^{\ast}$, \begin{equation}
Z(\eta^{\ast})=A+B(I-L_{\eta}^{\ast}D)^{-1}L_{\eta}^{\ast}C\label{transferfunction}\end{equation}
 where $L_{\eta}:H\rightarrow E^{\sigma}\otimes H$ is defined by
the formula $L_{\eta}h=\eta\otimes h$ (so $L_{\eta}^{\ast}(\theta\otimes h)=\tau(\langle\eta,\theta\rangle)h$). 
\end{itemize}
\end{theorem}

\begin{remark} Before giving the proof of Theorem \ref{realization},
we want to note that the result bears a strong resemblance to standard
results in the literature. We call special attention to \cite{AMcC99,AMcC02,jB00,BB04,BGM06a,BGM06b,BT98,BTV01}.
Indeed, we recommend \cite{jB00}, which is a survey that explains
the general strategy for proving the theorem. What is novel in our
approach is the adaptation of the results in the literature to accommodate
completely positive definite kernels.\end{remark}

Since the matrix in equation (\ref{transfermatrix}) and the function
in equation (\ref{transferfunction}) are familiar constructs in mathematical
systems theory, more particularly from $H^{\infty}$-control theory
(see, e.g., \cite{kZ}), we adopt the following terminology.

\begin{definition} Let $E$ be a $W^{\ast}$-correspondence over
a von Neumann algebra $M$. Suppose that $\sigma$ is a faithful normal
representation of $M$ on a Hilbert space $\mathcal{E}$ and that
$q_{1}$ and $q_{2}$ are projections in $\sigma(M)$. Then an operator
matrix $V=\left(\begin{array}{cc}
A & B\\
C & D\end{array}\right)$, where the entries $A$, $B$, $C$, and $D$, satisfy conditions
$(i)$ and $(ii)$ of Theorem \ref{realization} for some normal representation
$\tau$ of $\sigma(M)^{\prime}$ on a Hilbert space $H$, is called
a \emph{system matrix} provided $V$ is a coisometry (that is unitary,
if $E$ is full). If $V$ is a system matrix, then the function $A+B(I-L_{\eta}^{\ast}D)^{-1}L_{\eta}^{\ast}C$,
$\eta^{\ast}\in\mathbb{D}(E^{\sigma})^{\ast}$ is called the \emph{transfer
function determined by} $V$. \end{definition}

\begin{proof} As we just remarked, the hypothesis that $Z$ is a
Schur class function on $\Omega^{\ast}$ means that the kernel $k_{Z}$
in equation (\ref{cpkernel}) is completely positive definite in the
sense of \cite{BBLS}. Consequently, we may apply Theorem 3.2.3 of
\cite{BBLS}, which is a lovely extension of Kolmogorov's representation
theorem for positive definite kernels, to find an $N$-$B(\mathcal{E})$
$W^{\ast}$-correspondence $F$ and a function $\iota$ from $\Omega^{\ast}$
to $F$ such that $F$ is spanned by $N\iota(\Omega^{\ast})B(\mathcal{E})$
and such that for every $\eta_{1}$ and $\eta_{2}$ in $\Omega^{\ast}$
and every $a\in N$, \[
(id-Ad(Z(\eta_{1}^{\ast}),Z(\eta_{2}^{\ast})))\circ(id-\theta_{\eta_{1},\eta_{2}})^{-1}(a)=\langle\iota(\eta_{1}),a\iota(\eta_{2})\rangle.\]
 It follows that for every $b\in N$ and every $\eta_{1},\eta_{2}$
in $\Omega^{\ast}$, \[
b-Z(\eta_{1}^{\ast})bZ(\eta_{2}^{\ast})^{\ast}=\langle\iota(\eta_{1}),b\iota(\eta_{2})\rangle-\langle\iota(\eta_{1}),\langle\eta_{1},b\eta_{2}\rangle\iota(\eta_{2})\rangle\]
\[
=\langle\iota(\eta_{1}),b\iota(\eta_{2})\rangle-\langle\eta_{1}\otimes\iota(\eta_{1}),b\eta_{2}\otimes\iota(\eta_{2})\rangle.\]
 Thus, \begin{equation}
b+\langle\eta_{1}\otimes\iota(\eta_{1}),b\eta_{2}\otimes\iota(\eta_{2})\rangle=\langle\iota(\eta_{1}),b\iota(\eta_{2})\rangle+Z(\eta_{1}^{\ast})bZ(\eta_{2}^{\ast})^{\ast}.\label{eq}\end{equation}
 Set \[
G_{1}:=\overline{span}\{ bZ(\eta^{\ast})^{\ast}q_{2}T\oplus b\iota(\eta)q_{2}T\mid b\in N,\;\eta\in\Omega^{\ast},\; T\in B(\mathcal{E})\;\}\]
 and \[
G_{2}:=\overline{span}\{ bq_{2}T\oplus(b\eta\otimes\iota(\eta)q_{2}T)\mid b\in N,\;\eta\in\Omega^{\ast},\; T\in B(\mathcal{E})\;\}.\]
 Then $G_{1}$ is a sub $N$-$B(\mathcal{E})$ $W^{\ast}$-correspondence
of $B(\mathcal{E})\oplus F$ (where we use the assumption that $q_{2}Z(\eta^{\ast})=q_{2}Z(\eta^{\ast})q_{1}$)
and $G_{2}$ is a sub $N$-$B(\mathcal{E})$ $W^{\ast}$-correspondence
of $B(\mathcal{E})\oplus(E^{\sigma}\otimes_{N}F)$ . (The closure
in the definitions of $G_{1},G_{2}$ is in the $\sigma$-topology
of \cite{BDH88}. It then follows that $G_{1}$ and $G_{2}$ are $W^{\ast}$-correspondences
\cite[Consequences 1.8 (i)]{BDH88}). Define $v:G_{1}\rightarrow G_{2}$
by the equation \[
v(bZ(\eta^{\ast})^{\ast}q_{2}T\oplus b\iota(\eta)q_{2}T)=bq_{2}T\oplus(b\eta\otimes\iota(\eta)q_{2}T).\]
 It follows from (\ref{eq}) that $v$ is an isometry. It is also
clear that it is a bimodule map. We write $P_{i}$ for the orthogonal
projection onto $G_{i}$, $i=1,2$ and $\tilde{V}$ for the map \[
\tilde{V}:=P_{2}vP_{1}:q_{1}B(\mathcal{E})\oplus F\rightarrow q_{2}B(\mathcal{E})\oplus(E^{\sigma}\otimes_{N}F).\]
 Then $\tilde{V}$ is a partial isometry and, since $P_{1},v$ and
$P_{2}$ are all bimodule maps (see Lemma~\ref{bimodule}), so is
$\tilde{V}$. We write $\tilde{V}$ matricially: \[
\tilde{V}=\left(\begin{array}{cc}
\alpha & \beta\\
\gamma & \delta\end{array}\right),\]
 where $\alpha:q_{1}B(\mathcal{E})\rightarrow q_{2}B(\mathcal{E})$,
$\beta:F\rightarrow q_{2}B(\mathcal{E})$, $\gamma:q_{1}B(\mathcal{E})\rightarrow E^{\sigma}\otimes F$
and $\delta:F\rightarrow E^{\sigma}\otimes F$ and all these maps
are bimodule maps. Let $H_{0}$ be the Hilbert space $F\otimes_{B(\mathcal{E})}\mathcal{E}$
and note that $B(\mathcal{E})\otimes_{B(\mathcal{E})}\mathcal{E}$
is isomorphic to $\mathcal{E}$ (and the isomorphism preserves the
left $N$-action). Tensoring on the right by $\mathcal{E}$ (over
$B(\mathcal{E})$) we obtain a partial isometry \[
V_{0}:=\left(\begin{array}{cc}
A_{0} & B_{0}\\
C_{0} & D_{0}\end{array}\right):\left(\begin{array}{c}
\mathcal{E}_{1}\\
H_{0}\end{array}\right)\rightarrow\left(\begin{array}{c}
\mathcal{E}_{2}\\
E^{\sigma}\otimes H_{0}\end{array}\right).\]
 Here $A_{0}=\alpha\otimes I_{\mathcal{E}}$, $B_{0}=\beta\otimes I_{\mathcal{E}}$,
$C_{0}=\gamma\otimes I_{\mathcal{E}}$ and $D_{0}=\delta\otimes I_{\mathcal{E}}$.
These maps are well defined because the maps $\alpha,\beta,\gamma$
and $\delta$ are right $B(\mathcal{E})$-module maps. Since these
maps are also left $N$-module maps, so are $A_{0},B_{0},C_{0}$ and
$D_{0}$.

By the definition of $V_{0}$, its initial space is $G_{1}\otimes\mathcal{E}$
and its final space is $G_{2}\otimes\mathcal{E}$. In fact, $V_{0}$
induces an equivalence of the representations of $N$ on $G_{1}\otimes\mathcal{E}$
and on $G_{2}\otimes\mathcal{E}$.

It will be convenient to use the notation $K_{1}\preceq_{N}K_{2}$
if the Hilbert spaces $K_{1}$ and $K_{2}$ are both left $N$-modules
and the representation of $N$ on $K_{1}$ is equivalent to a subrepresentation
of the representation of $N$ on $K_{2}$. This means, of course,
that there is an isometry from $K_{1}$ into $K_{2}$ that intertwines
the two representations. If the two representations are equivalent
we write $K_{1}\simeq_{N}K_{2}$.

Using this notation, we can write $G_{1}\otimes\mathcal{E}\simeq_{N}G_{2}\otimes\mathcal{E}$.
Form $\mathcal{M}_{2}:=(\mathcal{E}_{2}\oplus(E^{\sigma}\otimes H_{0}))\ominus(G_{2}\otimes\mathcal{E})$,
which is a left $N$-module, and note that $L:=\mathcal{F}(E^{\sigma})\otimes\mathcal{M}_{2}$
also is a left $N$-module, where the representation of $N$ on $L$
is the induced representation. Since $L=\mathcal{F}(E^{\sigma})\otimes\mathcal{M}_{2}=\bigoplus_{n=0}^{\infty}((E^{\sigma})^{\otimes n}\otimes(\mathcal{M}_{2}))$,
it is evident that $(E^{\sigma}\otimes L)\oplus\mathcal{M}_{2}\simeq_{N}L$.
Indeed, the isomorphisms are just the natural ones that give the associativity
of the tensor products involved. Thus, $\mathcal{E}_{2}\oplus(E^{\sigma}\otimes(H_{0}\oplus L))=\mathcal{E}_{2}\oplus(E^{\sigma}\otimes H_{0})\oplus(E^{\sigma}\otimes L)=G_{2}\otimes\mathcal{E}\oplus\mathcal{M}_{2}\oplus E^{\sigma}\otimes L\simeq_{N}G_{2}\otimes\mathcal{E}\oplus L\simeq_{N}G_{1}\otimes\mathcal{E}\oplus L\preceq_{N}\mathcal{E}_{1}\oplus(H_{0}\oplus L)$.
Consequently, we obtain a coisometric operator $V:\mathcal{E}_{1}\oplus(H_{0}\oplus L)\rightarrow\mathcal{E}_{2}\oplus E^{\sigma}\otimes(H_{0}\oplus L)$
that intertwines the representations of $N$ and extends $V_{0}$.
Note that, if $V_{0}$ were known to be an isometry (so that $G_{2}\otimes\mathcal{E}\simeq_{N}G_{1}\otimes\mathcal{E}=\mathcal{E}_{1}\oplus H_{0}$
), then we would have equivalence above and $V$ can be chosen to
be unitary.

Assume that $E$ is full. We also write $\mathcal{M}_{1}$ for $(\mathcal{E}_{1}\oplus H_{0})\ominus G_{1}\otimes\mathcal{E}$.
Since $E$ is full, the representation $\rho$ of $N$ on $E^{\sigma}\otimes\mathcal{E}$
is faithful (Lemma~\ref{faithful}) and it follows that every representation
of $N$ is quasiequivalent to a subrepresentation of $\rho$. Write
$\mathcal{E}_{\infty}$ for the direct sum of infinitely many copies
of $\mathcal{E}$. Then $E^{\sigma}\otimes\mathcal{E}_{\infty}$ is
the direct sum of infinitely many copies of $E^{\sigma}\otimes\mathcal{E}$
and, thus, every representation of $N$ is equivalent to a subrepresentation
of the representation of $N$ on $E^{\sigma}\otimes\mathcal{E}_{\infty}$.
In particular, we can write $\mathcal{M}_{1}\oplus\mathcal{E}_{\infty}\preceq_{N}E^{\sigma}\otimes\mathcal{E}_{\infty}$.
Thus $\mathcal{E}_{1}\oplus(H_{0}\oplus\mathcal{E}_{\infty})=(G_{1}\otimes\mathcal{E})\oplus\mathcal{M}_{1}\oplus\mathcal{E}_{\infty}\preceq_{N}\mathcal{E}_{2}\oplus(E^{\sigma}\otimes H_{0})\oplus(E^{\sigma}\otimes\mathcal{E}_{\infty})=\mathcal{E}_{2}\oplus(E^{\sigma}\otimes(H_{0}\oplus\mathcal{E}_{\infty}))$.
So, replacing $H_{0}$ by $H_{0}\oplus\mathcal{E}_{\infty}$, we can
replace $V_{0}$ by an isometry and, using the argument just presented,
we conclude that the resulting $V$ is a unitary operator intertwining
the representations of $N$ and extending $V_{0}$.

So we let $V$ be the coisometry just constructed (and treat it as
unitary when $E$ is full). Writing $H:=H_{0}\oplus L$, we can express
$V$ in the matricial form as in part (iii) of the statement of the
theorem. Conditions (i) and (ii) then follow from the fact that $V$
intertwines the indicated representations of $N$. It is left to prove
(iv).

Setting $b=T=I$ in the definition of $v$ above and writing $v$
in a matricial form we see that \[
\left(\begin{array}{cc}
\alpha & \beta\\
\gamma & \delta\end{array}\right)\left(\begin{array}{c}
Z(\eta^{\ast})^{\ast}q_{2}\\
\iota(\eta)q_{2}\end{array}\right)=\left(\begin{array}{c}
q_{2}\\
\eta\otimes\iota(\eta)q_{2}\end{array}\right).\]
 Tensoring by $I_{\mathcal{E}}$ on the right and identifying $B(\mathcal{E})\otimes_{B(\mathcal{E})}\mathcal{E}$
with $\mathcal{E}$ as above, we find that \[
\left(\begin{array}{cc}
A_{0} & B_{0}\\
C_{0} & D_{0}\end{array}\right)\left(\begin{array}{c}
Z(\eta^{\ast})^{\ast}g\\
\iota(\eta)\otimes g\end{array}\right)=\left(\begin{array}{c}
g\\
\eta\otimes(\iota(\eta)\otimes g)\end{array}\right),\]
 for $g\in\mathcal{E}_{2}$. Since $A,B,C$ and $D$ extend $A_{0},B_{0},C_{0}$
and $D_{0}$ respectively, we can drop the subscript $0$. We also
use the fact that the matrix we obtain is a coisometry, and thus its
adjoint equals its inverse on its range. We conclude that \begin{equation}
\left(\begin{array}{cc}
A^{\ast} & C^{\ast}\\
B^{\ast} & D^{\ast}\end{array}\right)\left(\begin{array}{c}
g\\
\eta\otimes(\iota(\eta)\otimes g)\end{array}\right)=\left(\begin{array}{c}
Z(\eta^{\ast})^{\ast}g\\
\iota(\eta)\otimes g\end{array}\right).\label{eq3}\end{equation}
 Thus $\iota(\eta)\otimes g=B^{\ast}g+D^{\ast}(\eta\otimes(\iota(\eta)\otimes g))=B^{\ast}g+D^{\ast}L_{\eta}(\iota(\eta)\otimes g)$
and \[
\iota(\eta)\otimes g=(I-D^{\ast}L_{\eta})^{-1}B^{\ast}g.\]
 Combining this equality with the other equation that we get from
(\ref{eq3}), we have \[
Z(\eta^{\ast})^{\ast}g=A^{\ast}g+C^{\ast}L_{\eta}(I-D^{\ast}L_{\eta})^{-1}B^{\ast}g\;,\;\; g\in\mathcal{E}.\]
 Taking adjoints yields (iv). \end{proof}

Thus, Theorem \ref{realization} asserts that every Schur class function
determines a system matrix whose transfer function represents the
function. The system matrix is not unique in general, but as the proof
of Theorem \ref{realization} shows, it arises through a series of
natural choices. Of course, equation (\ref{transferfunction}) suggests
that every Schur class function represents an element in $H^{\infty}(E)$.
This is indeed the case, as the following converse shows.

\begin{theorem} \label{evaluationtheorem}Let $E$ be a $W^{\ast}$-correspondence
over a $W^{\ast}$-algebra $M$, and let $\sigma$ be a faithful normal
representation of $M$ on a Hilbert space $\mathcal{E}$. If $V=\left(\begin{array}{cc}
A & B\\
C & D\end{array}\right)$ is a system matrix determined by a normal representation $\tau$
of $N:=\sigma(M)^{\prime}$ on a Hilbert space $H$, then there is
an $X\in H^{\infty}(E)$, $\Vert X\Vert\leq1$, such that \[
\widehat{X}(\eta^{\ast})=A+B(I-L_{\eta}^{\ast}D)^{-1}L_{\eta}^{\ast}C\text{,}\]
 for all $\eta^{\ast}\in\mathbb{D}(E^{\sigma})^{\ast}$ and, conversely,
every $X\in H^{\infty}(E)$, $\Vert X\Vert\leq1$, may be represented
in this fashion for a suitable system matrix $V=\left(\begin{array}{cc}
A & B\\
C & D\end{array}\right)$. \end{theorem}

\begin{proof} For every $n\geq0$ we define an operator $K_{n}$
from $\mathcal{E}$ to $(E^{\sigma})^{\otimes n}\otimes\mathcal{E}$
as follows. For $n=0$, we set $K_{0}=A$ - an operator in $B(\mathcal{E})$.
For $n=1$, we define $K_{1}$, mapping $\mathcal{E}$ to $E^{\sigma}\otimes\mathcal{E}$,
to be $(I_{1}\otimes B)C$, where for all $k\geq1$, $I_{k}$ denotes
the identity operator on $(E^{\sigma})^{\otimes k}$. For $n\geq2$,
we set\[
K_{n}:=(I_{n}\otimes B)(I_{n-1}\otimes D)\cdots(I_{1}\otimes D)C.\]
 Note, first, that it follows from the properties of $A,B,C$ and
$D$ that, for every $n\geq0$ and every $a\in N$, $K_{n}a=(\varphi_{n}(a)\otimes I_{\mathcal{E}})K_{n}$
where $\varphi_{n}$ defines the left multiplication on $(E^{\sigma})^{\otimes n}$.
Thus, writing $\iota$ for the identity representation of $N$ on
$\mathcal{E}$, $K_{n}$ lies in the $\iota$-dual of $(E^{\sigma})^{\otimes n}$
which, by Theorem 3.6 and Lemma 3.7 of \cite{MS03}, is isomorphic
to $E^{\otimes n}$. Hence, for every $n\geq0$, $K_{n}$ defines
a unique element $\xi_{n}$ in $E^{\otimes n}$.

For every $n\geq0$ and $\eta\in E^{\sigma}$ we shall write $L_{n}(\eta)$
for the operator from $(E^{\sigma})^{\otimes n}\otimes\mathcal{E}$
to $(E^{\sigma})^{\otimes(n+1)}\otimes\mathcal{E}$ given by tensoring
on the left by $\eta$. Also note that, for $k\geq1$ and $n\geq0$,
$I_{k}\otimes K_{n}$ is an operator from $(E^{\sigma})^{\otimes k}\otimes\mathcal{E}$
to $(E^{\sigma})^{\otimes(k+n)}\otimes\mathcal{E}$. With this notation,
it is easy to see that, for all $k\geq1$ and $n\geq0$, \begin{equation}
(I_{k+1}\otimes K_{n})L_{k}(\eta)=L_{k+n}(\eta)(I_{k}\otimes K_{n}).\label{commute}\end{equation}
 Note, too, that we can write \[
\mathcal{F}(E^{\sigma})\otimes\mathcal{E}=\mathcal{E}\oplus(E^{\sigma}\otimes\mathcal{E})\oplus\cdots\oplus((E^{\sigma})^{\otimes m}\otimes\mathcal{E})\oplus\cdots\]
 and every operator on $\mathcal{F}(E^{\sigma})\otimes\mathcal{E}$
can be written in a matricial form with respect to this decomposition
(with indices starting at $0$). For every $m$, $0\leq m\leq\infty$,
we let $S_{m}$ be the operator defined by the matrix whose $i,j$
entry is $I_{j}\otimes K_{i-j}$, if $0\leq j\leq i\leq m$, and is
$0$ otherwise. (For $m=\infty$, it is not clear yet that the matrix
so constructed represents a bounded operator, but this will be verified
later).

So far we have not used the assumption that $V$ is a coisometry.
But if we take this into account, form the product $VV^{\ast}$, and
set it equal to $I_{\mathcal{E}\oplus(E^{\sigma}\otimes H)}$, we
find that \begin{align}
I_{\mathcal{E}}-AA^{\ast} & =BB^{\ast}\label{cond11}\\
CC^{\ast} & =I_{E^{\sigma}\otimes_{\tau}H}-DD^{\ast}\label{cond22}\\
AC^{\ast} & =-BD^{\ast}\label{cond12}\end{align}
 We claim that, for $1\leq j\leq i\leq m$, the following equations
hold, \begin{equation}
(I-S_{m}S_{m}^{\ast})_{i,j}=(I_{i}\otimes B)(I_{i-1}\otimes D)\cdots DD^{\ast}\cdots(I_{j-1}\otimes D^{\ast})(I_{j}\otimes B^{\ast});\label{comp1}\end{equation}
 that for $0<i\leq m$, \begin{equation}
(I-S_{m}S_{m}^{\ast})_{i,0}=(I_{i}\otimes B)(I_{i-1}\otimes D)\cdots DB^{\ast},\label{comp0}\end{equation}
 and that for $i=j=0$, \begin{equation}
(I-S_{m}S_{m}^{\ast})_{0,0}=BB^{\ast}.\label{comp00}\end{equation}
 Equation (\ref{comp00}) follows immediately from (\ref{cond11})
since $(S_{m})_{0,0}=A$. For $0<i\leq m$ we compute $(I-S_{m}S_{m}^{\ast})_{i,0}=-(S_{m})_{i,0}(S_{m})_{0,0}^{\ast}=-(I_{i}\otimes B)(I_{i-1}\otimes D)\cdots(I_{1}\otimes D)CA^{\ast}=(I_{i}\otimes B)(I_{i-1}\otimes D)\cdots(I_{1}\otimes D)DB^{\ast}$
where, in the last equality we used (\ref{cond12}). It is left to
prove (\ref{comp1}). Let us write $R_{i,j}$ for the left hand side
of (\ref{comp1}). (For $j=0<i$ we have $R_{i,0}=(I_{i}\otimes B)(I_{i-1}\otimes D)\cdots DB^{\ast}$
and when both are $0$, $R_{0,0}=BB^{\ast}$). We have $K_{0}K_{0}^{\ast}=AA^{\ast}=I-BB^{\ast}=I-R_{0,0}R_{0,0}^{\ast}$.
For $0=j<i\leq m$ we have $K_{i}K_{0}^{\ast}=(I_{i}\otimes B)(I_{i-1}\otimes D)\cdots(I_{1}\otimes D)CA^{\ast}=-(I_{i}\otimes B)(I_{i-1}\otimes D)\cdots(I_{1}\otimes D)DB^{\ast}=-R_{i,0}$
and for $0<j\leq i\leq m$, $K_{i}K_{j}^{\ast}=(I_{i}\otimes B)(I_{i-1}\otimes D)\cdots(I_{1}\otimes D)CC^{\ast}(I_{1}\otimes D^{\ast})\cdots(I_{j-1}\otimes D^{\ast})(I_{j}\otimes B^{\ast})=(I_{i}\otimes B)(I_{i-1}\otimes D)\cdots(I_{1}\otimes D)(I-DD^{\ast})(I_{1}\otimes D^{\ast})\cdots(I_{j-1}\otimes D^{\ast})(I_{j}\otimes B^{\ast})=(I_{i}\otimes B)(I_{i-1}\otimes D)\cdots(I_{1}\otimes D)(I_{1}\otimes D^{\ast})\cdots(I_{j-1}\otimes D^{\ast})(I_{j}\otimes B^{\ast})-(I_{i}\otimes B)(I_{i-1}\otimes D)\cdots(I_{1}\otimes D)DD^{\ast}(I_{1}\otimes D^{\ast})\cdots(I_{j-1}\otimes D^{\ast})(I_{j}\otimes B^{\ast})=I_{1}\otimes R_{i-1,j-1}-R_{i,j}$.

We have \[
(S_{m}S_{m}^{\ast})_{i,j}=\sum_{k=0}^{j}(S_{m})_{i,k}(S_{m})_{j,k}=\sum_{k=0}^{j}I_{k}\otimes K_{i-k}K_{j-k}^{\ast}=\sum_{l=0}^{j}I_{j-l}\otimes K_{i-j+l}K_{l}^{\ast}.\]
 Using the computation above, we get, for $i=j\leq m$, \[
(S_{m}S_{m}^{\ast})_{i,i}=I_{i}\otimes(I-R_{0,0}R_{0,0}^{\ast})+\sum_{l=1}^{i}(I_{i-l+1}\otimes R_{l-1,l-1}-I_{i-l}\otimes R_{l,l})=I-R_{i,i}\]
 and, for $j<i\leq m$, \[
(S_{m}S_{m}^{\ast})_{i,j}=-I_{j}\otimes R_{i-j,0}+\sum_{l=1}^{j}(I_{j-l+1}\otimes R_{i-j+l-1,l-1}-I_{j-l}\otimes R_{i-j+l,l})=-R_{i,j}.\]
 This completes the proof of the claim. If we let $R$ be the operator
whose matrix is $(R_{i,j})$ (letting $R_{i,j}=0$ if $i$ or $j$
is larger than $m$) then we get $R=I-S_{m}S_{m}^{\ast}$. But it
is easy to verify that $R$ is a positive operator and, thus, $\Vert S_{m}\Vert\leq1$.
This holds for every $m$ and, therefore, we can find a weak limit
point of the sequence $\{ S_{m}\}$. But this limit point it clearly
equal to $S_{\infty}$, showing that $S_{\infty}$ is indeed a bounded
operator, with norm at most $1$.

Recall that the induced representation of $H^{\infty}(E)$ on $\mathcal{F}(E)\otimes_{\sigma}\mathcal{E}$
is the representation that maps $X\in H^{\infty}(E)$ to $\sigma^{\mathcal{F}(E)}(X):=X\otimes I_{\mathcal{E}}$.
The representation is faithful and is a homeomorphism with respect
to the ultraweak topologies. Its image is the ultraweakly closed subalgebra
of $B(\mathcal{F}(E)\otimes\mathcal{E})$ generated by the operators
$T_{\xi}\otimes I_{\mathcal{E}}$ and $\varphi_{\infty}(a)\otimes I_{\mathcal{E}}$
for $\xi\in E$ and $a\in M$. Similarly one defines the induced representation
$\iota^{\mathcal{F}(E^{\sigma})}$ of $H^{\infty}(E^{\sigma})$ on
$\mathcal{F}(E^{\sigma})\otimes\mathcal{E}$ and its image is generated
by the operators $T_{\eta}\otimes I$ and $\varphi_{\infty}(b)\otimes I$
for $\eta\in E^{\sigma}$ and $b\in N$. Recall also, from \cite[Theorem 3.9]{MS03},
that there is a unitary operator $U:\mathcal{F}(E^{\sigma})\otimes\mathcal{E}\rightarrow\mathcal{F}(E)\otimes\mathcal{E}$
such that \[
(\iota^{\mathcal{F}(E^{\sigma})}(H^{\infty}(E^{\sigma})))^{\prime}=U^{\ast}\sigma^{\mathcal{F}(E)}(H^{\infty}(E))U.\]

That is, $U$ gives an explicit representation of $H^{\infty}(E^{\sigma})$
as the commutant of the induced algebra $\sigma^{\mathcal{F}(E)}(H^{\infty}(E))$.
Thus, to show that $S_{\infty}=U^{\ast}(X\otimes I)U$ for an $X\in H^{\infty}(E)$,
we need only show that $S_{\infty}$ lies in the commutant of $\iota^{\mathcal{F}(E^{\sigma})}(H^{\infty}(E^{\sigma}))$.
And for this, we only have to show that it commutes with the operators
$\varphi_{\infty}(b)\otimes I$, $b\in N$, and $T_{\eta}\otimes I$,
$\eta\in E^{\sigma}$. Note that, matricially, $\varphi_{\infty}(b)\otimes I$
is a diagonal operator whose $i,i$ entry is $\varphi_{i}(b)$. For
$S_{\infty}$ to commute with it we should have, for all $j\leq i$,
\[
(I_{j}\otimes K_{i-j})(\varphi_{j}(b)\otimes I)=(\varphi_{i}(b)\otimes I)(I_{j}\otimes K_{i-j}).\]
 This equality is obvious for $j>0$. For $j=0$ it amounts to the
equality \[
K_{i}b=(\varphi_{i}(b)\otimes I_{\mathcal{E}})K_{i}\]
 and, this, as was mentioned above, follows immediately from the properties
of $A,B,C$ and $D$. To show that $S_{\infty}$ commutes with every
$T_{\eta}\otimes I$, $\eta\in E^{\sigma}$, note that, matricially,
the $i,j$ entry of $T_{\eta}\otimes I$ vanishes unless $i=j+1$
and, in this case the entry is $L_{j}(\eta)$. Equation (\ref{commute})
then ensures that $S_{\infty}$ and $T_{\eta}\otimes I$ commute.

Thus, by \cite[Theorem 3.9]{MS03}, there is an element $X\in H^{\infty}(E)$
such that $S_{\infty}=U^{\ast}(X\otimes I)U$ ($=U^{\ast}\sigma^{\mathcal{F}(E)}(X)U$).
Since $S_{\infty}$ has norm at most one, so does $X$.

It remains to show that $X$ is given by the transfer function built
from $V$. To this end, fix $\xi\in E$ and recall that $\xi$ defines
a map $W(\xi):\mathcal{E}\rightarrow E^{\sigma}\otimes\mathcal{E}$
via the formula $W(\xi)^{\ast}(\eta\otimes h)=L_{\xi}^{\ast}\eta(h)$,
$\eta\otimes h\in E^{\sigma}\otimes\mathcal{E}$ (See \cite[Theorem 3.6]{MS03}.),
and that $W$ maps $E$ onto the $\iota$-dual of $E^{\sigma}$. The
desired properties follow easily from the definition of $W$. For
every $k\geq0$, $I_{k}\otimes W(\xi)^{\ast}$ is a map from $(E^{\sigma})^{\otimes k+1}\otimes\mathcal{E}$
into $(E^{\sigma})^{\otimes k}\otimes\mathcal{E}$. An easy computation
shows that it is equal to the restriction of $U^{\ast}(T_{\xi}^{\ast}\otimes I_{\mathcal{E}})U$
to $(E^{\sigma})^{\otimes k+1}\otimes\mathcal{E}$. (Recall from \cite[Lemma 3.8]{MS03}
that the restriction of $U$ to $(E^{\sigma})^{\otimes k+1}\otimes\mathcal{E}$
is defined by the equation $U(\eta_{1}\otimes\cdots\otimes\eta_{k+1}\otimes h)=(I_{k}\otimes\eta_{1})\cdots(I_{1}\otimes\eta_{k})\eta_{k+1}(h)$.)

It then follows that the $i,j$ entry of the matrix associated with
$U^{\ast}(T_{\xi}\otimes I_{\mathcal{E}})U$ vanishes unless $i=j+1$
and \[
(U^{\ast}(T_{\xi}\otimes I_{\mathcal{E}})U)_{j+1,j}=I_{j}\otimes W(\xi).\]
 Similarly one can show that, for $\xi\in E^{\otimes k}$, the $i,j$
entry of the matrix associated with $U^{\ast}(T_{\xi}\otimes I_{\mathcal{E}})U$
vanishes unless $i=j+k$ and \[
(U^{\ast}(T_{\xi}\otimes I_{\mathcal{E}})U)_{j+k,j}=I_{j}\otimes W(\xi).\]
 In the last equation, $W(\xi)$, $\xi\in E^{\otimes k}$, is a map
from $\mathcal{E}$ to $(E^{\sigma})^{\otimes k}\otimes\mathcal{E}$.

Recall that we defined $\xi_{n}$ to be the vectors in $E^{\otimes n}$
with $W(\xi_{n})=K_{n}$. Thus we see that the $n^{th}$ lower diagonal
in the matricial form of $S_{\infty}$ is the matricial form of $U^{\ast}(T_{\xi_{n}}\otimes I_{\mathcal{E}})U$.

Recall from the discussion at the end of Section 2 in \cite{MS03}
that $S_{\infty}$ is the ultraweak limit of the sequence $\Sigma_{k}$
where \[
\Sigma_{k}=\sum_{j=0}^{k-1}(1-\frac{j}{k})U^{\ast}(T_{\xi_{j}}\otimes I)U.\]
 Hence $X$ is the ultraweak limit of $X_{k}$ where \[
X_{k}=\sum_{j=0}^{k-1}(1-\frac{j}{k})T_{\xi_{j}}\]
 and, for $\eta\in E^{\sigma}$, $\widehat{X}(\eta^{\ast})$ is the
ultraweak limit of $\widehat{X}_{k}(\eta^{\ast})=\sum_{j=0}^{k-1}(1-\frac{j}{k})\widehat{T_{\xi_{j}}}(\eta^{\ast})$.

Fix $\eta\in E^{\sigma}$ and $k\geq1$. Then it is easy to check
that, in the notation of the theorem, $L_{\eta}^{\ast}(I_{k}\otimes B)=(I_{k-1}\otimes B)L_{\eta}^{\ast}$
and $L_{\eta}^{\ast}(I_{k}\otimes D)=(I_{k-1}\otimes D)L_{\eta}^{\ast}$,
all as operators on $(E^{\sigma})^{\otimes k}\otimes H$. It then
follows that for $n\geq1$, \[
(L_{\eta}^{\ast})^{n}W(\xi_{n})=(L_{\eta}^{\ast})^{n}K_{n}=B(L_{\eta}^{\ast}D)^{n-1}L_{\eta}^{\ast}C\]
 and \[
A+B(I-L_{\eta}^{\ast}D)^{-1}L_{\eta}^{\ast}C=A+\sum_{n=1}^{\infty}B(L_{\eta}^{\ast}D)^{n-1}L_{\eta}^{\ast}C=\sum_{n=0}^{\infty}(L_{\eta}^{\ast})^{n}W(\xi_{n}).\]
 (Note that the last series converges in norm). It follows from \cite[Proposition 5.1]{MS03}
that $\widehat{T_{\xi_{n}}}(\eta^{\ast})=(L_{\eta}^{\ast})^{n}W(\xi_{n})$
and, thus, we finally conclude that $\widehat{X}(\eta^{\ast})=A+B(I-L_{\eta}^{\ast}D)^{-1}L_{\eta}^{\ast}C$.

The `converse' portion of the Theorem is immediate from Theorems \ref{SchurHinfty}
and \ref{realization}. \end{proof}

\begin{corollary} \label{extension} Every Schur class operator function
defined on a subset $\Omega^{\ast}$ of $\mathbb{D}(E^{\sigma})^{\ast}$
with values in some $B(\mathcal{E})$ can be extended to a Schur class
operator function defined on all of $\mathbb{D}(E^{\sigma})^{\ast}$.
\end{corollary}

\begin{proof} Let $Z$ be a Schur class function on $\Omega^{\ast}$
and apply Theorem~\ref{realization} to represent $Z$ as the restriction
to $\Omega^{\ast}$ of a transfer function. The result then follows
from the evident combination of Theorems \ref{evaluationtheorem}
and \ref{SchurHinfty}. \end{proof}

Recall that every element $X$ in $H^{\infty}(E)$ with $\Vert X\Vert\leq1$
defines a Schur class operator function by evaluation at $\eta^{\ast}$
for $\eta\in\mathbb{D}(E^{\sigma})$ (where $\sigma$ is a suitable
prescribed faithful normal representation of $M$) . We usually suppress
reference to $\sigma$ and write $\widehat{X}$ for this Schur class
operator function. In general, however, the map $X\rightarrow\widehat{X}$
is not one-to-one, and whether it is or not depends on the choice
of $\sigma$. Indeed, in the particular case when $M=\mathbb{C}$
and $E=\mathbb{C}^{n}$, so $H^{\infty}(E)$ is $\mathcal{L}_{n}$,
and when $\sigma$ is the identity representation of $\mathbb{C}$,
Davidson and Pitts showed that the kernel of the map $X\mapsto\widehat{X}$
is precisely the commutator ideal in $\mathcal{L}_{n}$ \cite{DP98}.
We shall show in the next lemma that given $E$, if $\sigma$ is chosen
to be faithful and have \emph{infinite uniform multiplicity}, meaning
that $\sigma$ is an infinite multiple of another faithful normal
representation of $M$, then the map $X\mapsto\widehat{X}$ will be
one-to-one. It will be convenient to write $K(\sigma)$ for the kernel
of the map determined by $\sigma$, so that \begin{align}
K(\sigma) & =\{ X\in H^{\infty}(E):\widehat{X}(\eta^{\ast})=0,\;\;\eta\in\mathbb{D}(E^{\sigma})\}\label{KernalSigma}\\
 & =\{ X\in H^{\infty}(E):\sigma\times\eta^{\ast}(X)=0,\;\;\eta\in\mathbb{D}(E^{\sigma})\}.\nonumber \end{align}

\begin{lemma} \label{1to1}If $\sigma$ is a faithful normal representation
of $M$ on a Hilbert space $H$ of infinite multiplicity, then $K(\sigma)=0$.
Moreover, if $\{ X_{\beta}\}$ is a bounded net in $H^{\infty}(E)$
and if there is an element $X\in H^{\infty}(E)$ such that for every
$\eta\in\mathbb{D}(E^{\sigma})$, $\widehat{X}_{\beta}(\eta^{\ast})\rightarrow\widehat{X}(\eta^{\ast})$
in the weak operator topology, then $X_{\beta}\rightarrow X$ ultraweakly.
\end{lemma}

\begin{proof} It follows from the structure of isomorphisms of von
Neumann algebras that any two infinite multiples of faithful representations
of a von Neumann algebra are unitarily equivalent. It follows, therefore,
that to prove the lemma, we can pick a special representation with
this property that is convenient for our purposes. So let $\pi$ be
the representation of $M$ on $\mathcal{F}(E)\otimes_{\sigma}H$ defined
by $\pi=\varphi_{\infty}\otimes I_{H}$. We shall see that $K(\pi)=\{0\}$.
For $\xi\in E$ let $V(\xi)=T_{\xi}\otimes I_{H}$. Then $(V,\pi)$
is a representation of $E$ on $\mathcal{F}(E)\otimes_{\sigma}H$.
The integrated form of this representation is the induced representation
$\pi^{\mathcal{F}(E)}$ restricted to $H^{\infty}(E)$. It is a faithful
representation of $H^{\infty}(E)$. For $0\leq r\leq1$, $(rV,\pi)$
is also a representation of $E$. It follows from \cite[Lemma 7.11]{MS03}
that, for every $X\in H^{\infty}(E)$, the limit in the strong operator
topology of $(\pi\times rV)(X)$, as $r\rightarrow1$, is $(\pi\times V)(X)$.
Thus, for $X\neq0$ in $H^{\infty}(E)$, there is an $r$, $0\leq r<1$,
such that $(\pi\times rV)(X)\neq0$. Since for such $r$ the inequality
$\Vert rV\Vert<1$ holds, and we conclude that $K(\pi)=\{0\}$.

For the second assertion of the lemma, suppose a bounded net $\{ X_{\beta}\}$
in $H^{\infty}(E)$ has the property that for every $\eta\in\mathbb{D}(E^{\pi})$,
$\widehat{X}_{\beta}(\eta^{\ast})\rightarrow0$. Since the net is
bounded, it has a ultraweak limit point $X_{0}$ in $H^{\infty}(E)$.
Since {}``evaluation at $\eta^{\ast}$\textquotedblright\ is the
same as applying a ultraweakly continuous representation, we see that
$\widehat{X}_{\beta}(\eta^{\ast})\rightarrow\widehat{X}_{0}(\eta^{\ast})$
for every $\eta\in\mathbb{D}(E^{\pi})$. But then, $\widehat{X}_{0}(\eta^{\ast})=0$
for every $\eta\in\mathbb{D}(E^{\pi})$ and, consequently, $X_{0}=0$
by the first assertion of the lemma. \end{proof}

With this lemma in hand, we summarize the results of this section
for future reference in the following corollary.

\begin{corollary} \label{SchurMult}Let $E$ be a $W^{\ast}$-correspondence
over the $W^{\ast}$-algebra $M$, let $\sigma$ be a faithful normal
representation of $M$ on the Hilbert space $\mathcal{E}$ and assume
that $\sigma$ has infinite multiplicity. Then the map $X\rightarrow\widehat{X}$
is a bijection from the closed unit ball of $H^{\infty}(E)$ onto
the space of Schur class $B(\mathcal{E})$-valued functions on $\mathbb{D}(E^{\sigma})^{\ast}$.
Further, for each $X$ in the closed unit ball of $H^{\infty}(E)$,
$\widehat{X}$ is the transfer function associated with a system matrix
$V=\left(\begin{array}{cc}
A & B\\
C & D\end{array}\right)$ defined in terms of a suitable auxiliary normal representation $\tau$
of $\sigma(M)^{\prime}$ on a Hilbert space $H$, and conversely,
each such transfer function on $\mathbb{D}(E^{\sigma})^{\ast}$, \[
\eta^{\ast}\rightarrow A+B(I-L_{\eta}^{\ast}D)^{-1}L_{\eta}^{\ast}C\text{,}\]
 is of the form $\widehat{X}$ for a uniquely determined $X\in H^{\infty}(E)$:
$\widehat{X}(\eta^{\ast})=A+B(I-L_{\eta}^{\ast}D)^{-1}L_{\eta}^{\ast}C$
for all $\eta\in\mathbb{D}(E^{\sigma})$. \end{corollary}

\begin{proof} The proof is just the evident combination of Lemma
\ref{1to1} and Theorems \ref{SchurHinfty}, \ref{realization}, and
\ref{evaluationtheorem}. \end{proof}

\begin{remark} One may well wonder why not stipulate at the outset
that all $\sigma$'s have uniform infinite multiplicity. It turns
out that in many interesting examples, such as those coming from graphs,
which we discuss in the last section, the principal $\sigma$'s one
wants to consider fail to have this property. \end{remark}

\section{Applications to automorphisms of the Hardy algebra}

In this section we apply the analysis of Schur class functions to
study automorphisms of $H^{\infty}(E)$. Our first goal is to show
that under very general assumptions, the automorphisms are obtained
by composition with (certain) biholomorphic automorphisms of the open
unit ball of the dual correspondence. For the case were $E=\mathbb{C}^{n}\,$,
so that $H^{\infty}(E)$ is the algebra $\mathcal{L}_{n}$ studied
by Davidson and Pitts and by Popescu, this was shown for the dual
correspondence associated with the one dimensional representation
$\sigma$ of $\mathbb{C}$ by Davidson and Pitts in \cite{DP98}.

Throughout this section we will focus on automorphisms $\alpha$ of
$H^{\infty}(E)$ that are \emph{completely isometric} and $w^{\ast}$\emph{-homeomorphisms}.
Also, we shall usually assume that the restriction of $\alpha$ to
$\varphi_{\infty}(M)$ is the identity.

It is known that, in various settings, one can assume much less. In
\cite{DP98}, the authors begin by assuming that $\alpha$ is simply
an algebraic automorphism but, to get the one-to-one correspondence
with automorphisms of the unit ball of the dual, they need to impose
also the assumption that the automorphism is contractive. It then
follows from their results that it is, in fact, completely isometric
and a $w^{\ast}$-homeomorphism. In \cite{KK04}, Katsoulis and Kribs
show that in the setting when $E$ is determined by a directed graph,
$G$ say, so $H^{\infty}(E)$ is the algebra they denote by $\mathcal{L}_{G}$,
an algebraic automorphism is always norm-continuous and $w^{\ast}$-continuous.

As for the assumption that the restriction of $\alpha$ to $\varphi_{\infty}(M)$
is the identity, we shall see that for many purposes this is no significant
restriction. However, in some situations, it can be a significant
technical headache to sort out what happens if we don't impose the
assumption. We will comment on this further, as we proceed. (See,
in particular, Remark~\ref{nonid}).

So, for the remainder of this section, unless specified otherwise,
$E$ will be a fixed $W^{\ast}$-correspondence over a $W^{\ast}$-algebra
$M$ and $\alpha$ will be a fixed automorphism of $H^{\infty}(E)$
that is completely isometric, $w^{\ast}$-homeomorphic and fixes $\varphi_{\infty}(M)$
element-wise. Also, $\sigma$ will be a faithful normal $\ast$-representation
of $M$ on a Hilbert space $H$.

We think about elements of $H^{\infty}(E)$ as functions on $\mathbb{D}(E^{\sigma})^{\ast}$
via the functional representation developed in the preceding section
and we want to study the transposed action of $\alpha$ on $\mathbb{D}(E^{\sigma})^{\ast}$.
For every $\eta\in\mathbb{D}(E^{\sigma})$, let $\tau(\eta):H\rightarrow E\otimes_{\sigma}H$
be defined by the equation \begin{equation}
\tau(\eta)^{\ast}(\xi\otimes h)=\widehat{\alpha(T_{\xi})}(\eta^{\ast})h\;(=(\sigma\times\eta^{\ast})(\alpha(T_{\xi}))h\;)\text{,}\label{tau}\end{equation}
 $\xi\otimes h\in E\otimes_{\sigma}H$. (Observe that if $\alpha$
is the identity automorphism of $H^{\infty}(E)$, then this equation
implies that $\tau$ is the identity map, as it should.) The next
lemma shows that $\tau(\eta)$ is well defined and is an element in
the closed unit ball of $E^{\sigma}$. Thus $\tau$ is a map from
$\mathbb{D}(E^{\sigma})$ into $\overline{\mathbb{D}(E^{\sigma})}$.
What we would really like to show, however, is that $\tau$ carries
$\mathbb{D}(E^{\sigma})$ into $\mathbb{D}(E^{\sigma})$, not the
closure. At this stage, we can only arrange for this under special
circumstances: Theorem \ref{strictcontraction} below. The restriction
on circumstances, however, is not so limiting as to eliminate many
interesting examples. We also want to show that $\tau$ is holomorphic
on $\mathbb{D}(E^{\sigma})$ in the usual sense of infinite dimensional
holomorphy \cite{HP57}.

\begin{lemma} \label{welldef}For each $\eta\in\mathbb{D}(E^{\sigma})$,
$\tau(\eta)$ is well defined and lies in the closed unit ball of
$E^{\sigma}$. \end{lemma}

\begin{proof} For $\xi\in E$, let $S(\xi):=(\sigma\times\eta^{\ast})(\alpha(T_{\xi}))$.
For every $a,b\in M$, $S(a\xi b)=(\sigma\times\eta^{\ast})(\alpha(T_{a\xi b}))=(\sigma\times\eta^{\ast})(\alpha(\varphi_{\infty}(a)T_{\xi}\varphi_{\infty}(b)))=(\sigma\circ\alpha)(\varphi_{\infty}(a))(\sigma\times\eta^{\ast})(\alpha(T_{\xi}))(\sigma\circ\alpha)(\varphi_{\infty}(b))$.
By our assumption, $\sigma\circ\alpha\circ\varphi_{\infty}=\sigma\circ\varphi_{\infty}$
and, thus, $(S,\sigma)$ is a covariant pair. Also, $S$ is a completely
contractive map of $E$ into $B(H)$ as a composition of three completely
contractive maps. Thus $\tilde{S}^{\ast}=\tau(\eta)$ lies in the
closed unit ball of $E^{\sigma}$. \end{proof}

To determine circumstances under which $\tau$ maps $\mathbb{D}(E^{\sigma})$
into $\mathbb{D}(E^{\sigma})$, we fix $\eta\in\mathbb{D}(E^{\sigma})$
and determine circumstances under which $\tau(z\eta)\in\mathbb{D}(E^{\sigma})$,
for every $z\in\mathbb{D}$ :$=\{ z\in\mathbb{C}\mid|z|<1\}$. This
will prove that $\tau$ maps $\mathbb{D}(E^{\sigma})$ into itself.

So for $z\in\mathbb{D}$, we define \begin{equation}
F(z):=\tau(\bar{z}\eta)^{\ast}.\label{F}\end{equation}
 Thus, $F(z)(\xi\otimes h)=(\sigma\times z\eta^{\ast})(\alpha(T_{\xi}))h$
for $\xi\in E$ and $h\in H$.

\begin{lemma} \label{analytic}$F$ is an analytic function from
$\mathbb{D}$ into $B(E\otimes H,H)$. \end{lemma}

\begin{proof} Fix $\xi\otimes h\in E\otimes H$ with $\Vert\xi\Vert\leq1$
and $k\in H$, and consider the expression \[
\langle F(z)(\xi\otimes h),k\rangle=\langle\widehat{\alpha(T_{\xi})}(z\eta^{\ast})h,k\rangle.\]
 Since $\alpha(T_{\xi})\in H^{\infty}(E)$ and $\Vert\alpha(T_{\xi})\Vert\leq1$,
we know from Theorem~\ref{evaluationtheorem} that we can write $\widehat{\alpha(T_{\xi})}(z\eta^{\ast})=A+B(I-zL_{\eta}^{\ast}D)^{-1}zL_{\eta}^{\ast}C$
for some system matrix. Thus \[
\widehat{\alpha(T_{\xi})}(z\eta^{\ast})=A+zBL_{\eta}^{\ast}C+\sum_{k=2}^{\infty}z^{k}B(L_{\eta}^{\ast})^{k-1}L_{\eta}^{\ast}C.\]
 Hence, for every $\xi\otimes h\in E\otimes H$ (even when $\Vert\xi\Vert>1$)
and $k\in H$, the function $z\mapsto\langle F(z)(\xi\otimes h),k\rangle$
is analytic. Since $\Vert F(z)\Vert\leq1$ by Lemma~\ref{welldef},
$|\langle F(z)g,k\rangle|\leq\Vert g\Vert\Vert k\Vert$ for every
$g\in E\otimes H$ and $k\in H$ and it follows that, for each such
$g,k$, the function $f_{g,k}(z):=\langle F(z)g,k\rangle$ is analytic
in $\mathbb{D}$ and $|f_{g,k}(z)|\leq\Vert g\Vert\Vert k\Vert$.
We can then write $f_{g,k}$ as a convergent power series $f_{g,k}(z)=\sum_{k=0}^{\infty}a_{n}(g,k)z^{n}$
and, for every $n\geq0$, $|a_{n}(g,k)|\leq\Vert g\Vert\Vert k\Vert$.
But then there are operators $A_{n}\in B(E\otimes H,H)$ with $\Vert A_{n}\Vert\leq1$
such that $a_{n}(g,k)=\langle A_{n}g,k\rangle$ for $g\in E\otimes H$
and $k\in H$. Hence $F(z)=\sum_{k=0}^{\infty}z^{n}A_{n}$ where the
sum converges in the weak operator topology. Since $|z|<1$ and the
norms of $\{ A_{n}\}$ are bounded by $1$, the series converges to
$F(z)$, for $z\in\mathbb{D}$, in the norm topology. We conclude
that $F(z)$ is analytic. \end{proof}

If we were dealing with scalar-valued functions, we would be able
to assert that $|F(z)|<1$ for all $z\in\mathbb{D}$, unless $F$
is constant, by the maximum modulus theorem. Unfortunately, an unalloyed
version of the maximum modulus theorem does not hold in our setting.
\ This is what leads to the special hypotheses on $\tau$ in Theorem
\ref{strictcontraction}. The next few results, then, which lead up
to Theorem \ref{strictcontraction} come out of our efforts to find
a serviceable replacement for the maximum modulus theorem. Our first
theorem in this direction, Theorem \ref{decomp}, is closely related
to \cite[Proposition V.2.1]{SzNF70}. It does not seem to follow directly
from this result, however. Instead, we appeal to the following lemma,
which in turn is an immediate application of an operator form of the
classical Pick criterion for interpolating operators at pre-assigned
points by operator-valued analytic functions. As such, it may be traced
back to Sz.-Nagy and Koranyi's influential paper \cite{SzNK56}. It
also is a consequence of Theorem 6.2 in \cite{MS03}, where it is
presented as a corollary of our Nevanlinna-Pick Theorem.

\begin{lemma} \label{matrix} If $K,H$ are Hilbert spaces and if
$F:\mathbb{D}\rightarrow B(K,H)$ is an analytic function satisfying
$\Vert F(z)\Vert\leq1$ for all $z\in\mathbb{D}$, then, for every
$z_{1},z_{2}\in\mathbb{D}$, the matrix \[
\left(\begin{array}{cc}
\frac{I_{H}-F(z_{1})F(z_{1})^{\ast}}{1-|z_{1}|^{2}} & \frac{I_{H}-F(z_{1})F(z_{2})^{\ast}}{1-z_{1}\bar{z_{2}}}\\
\frac{I_{H}-F(z_{2})F(z_{1})^{\ast}}{1-z_{2}\bar{z_{1}}} & \frac{I_{H}-F(z_{2})F(z_{2})^{\ast}}{1-|z_{2}|^{2}}\end{array}\right)\]
 is positive. In particular (setting $z_{1}=z$ and $z_{2}=0$), for
every $z\in\mathbb{D}$, \begin{equation}
\left(\begin{array}{cc}
\frac{I_{H}-F(z)F(z)^{\ast}}{1-|z|^{2}} & I_{H}-F(z)F(0)^{\ast}\\
I_{H}-F(0)F(z)^{\ast} & I_{H}-F(0)F(0)^{\ast}\end{array}\right)\geq0.\label{2x2}\end{equation}

\end{lemma}

\begin{theorem} \label{decomp} Suppose $H$ and $K$ are Hilbert
spaces and suppose $F:\mathbb{D}\rightarrow B(K,H)$ is an analytic
function that satisfies the following conditions:

\begin{enumerate}
\item [(1)] $\Vert F(z)\Vert\leq1$ for all $z\in\mathbb{D}$.
\item [(2)] There are projections $P_{1},P_{2}$ in $B(H)$ that sum to
$I_{H}$ and projections $Q_{1},Q_{2}$ in $B(K)$ that sum to $I_{K}$
and satisfy:

\begin{enumerate}
\item [(i)] $P_{1}F(0)Q_{2}=0$ and $P_{2}F(0)Q_{1}=0$.
\item [(ii)] $P_{1}F(0)F(0)^{*}P_{1}=P_{1}$.
\item [(iii)] $P_{2}F(0)F(0)^{*}P_{2}\leq rP_{2}$ for some $0<r<1$. 
\end{enumerate}
\end{enumerate}
Then, for every $z\in\mathbb{D}$,

\begin{enumerate}
\item [(1)] $P_{1}F(z)Q_{2}=0$.
\item [(2)] $P_{1}F(z)Q_{1}=P_{1}F(0)Q_{1}$.
\item [(3)] There is a function $q_{0}(z)$ on $\mathbb{D}$, such that
$0<q_{0}(z)<1$ for all $z\in\mathbb{D}$, and such that $P_{2}F(z)F(z)^{\ast}P_{2}\leq q_{0}(z)P_{2}$. 
\end{enumerate}
\end{theorem}

\begin{proof} It will be convenient to use the projections $P_{1},P_{2}$
and $Q_{1},Q_{2}$ to write $F(z)$ matricially as \[
F(z)=\left(\begin{array}{cc}
A(z) & B(z)\\
C(z) & D(z)\end{array}\right)\]
 so that, by assumption, \[
F(0)=\left(\begin{array}{cc}
A(0) & 0\\
0 & D(0)\end{array}\right)\]
 where $A(0)A(0)^{*}=P_{1}$ and $D(0)D(0)^{*}\leq rP_{2}$.

Since $F$ satisfies the conditions of Lemma~\ref{matrix}, Equation~\ref{2x2}
holds for all $z\in\mathbb{D}$. Compressing each entry of the matrix
in (\ref{2x2}) to the range of $P_{1}$ and using the fact that $A(0)A(0)^{*}=P_{1}$
and that $P_{1}F(0)Q_{2}=0$, we get \begin{equation}
\left(\begin{array}{cc}
\frac{P_{1}-P_{1}F(z)F(z)^{*}P_{1}}{1-|z|^{2}} & P_{1}-P_{1}F(z)Q_{1}A(0)^{*}\\
P_{1}-A(0)Q_{1}F(z)^{*}P_{1} & 0\end{array}\right)\geq0.\end{equation}
 It follows that $P_{1}=P_{1}F(z)Q_{1}A(0)^{*}$. Thus $0\leq(P_{1}F(z)Q_{1}-A(0))(Q_{1}F(z)^{*}P_{1}-A(0)^{*})=P_{1}F(z)Q_{1}F(z)^{*}P_{1}+A(0)A(0)^{*}-P_{1}F(z)Q_{1}A(0)^{*}-A(0)Q_{1}F(z)^{*}P_{1}\leq0$.
Consequently, $A(0)=P_{1}F(z)Q_{1}$ (for every $z\in\mathbb{D}$).

But then $P_{1}F(z)Q_{1}F(z)^{*}P_{1}=P_{1}$ and, since $P_{1}F(z)F(z)^{*}P_{1}\leq\newline P_{1}$,
$P_{1}F(z)Q_{2}=0$. This proves (1) and (2).

Compress each entry of (\ref{2x2}) to the range of $P_{2}$ to get
\begin{equation}
\left(\begin{array}{cc}
\frac{P_{2}-P_{2}F(z)F(z)^{\ast}P_{2}}{1-|z|^{2}} & P_{2}-P_{2}F(z)Q_{2}D(0)^{\ast}\\
P_{2}-D(0)Q_{2}F(z)^{\ast}P_{2} & P_{2}-D(0)D(0)^{\ast}\end{array}\right)\geq0.\label{pos}\end{equation}
 Write $\Delta$ for the positive square root of $P_{2}-D(0)D(0)^{\ast}$
and note that $\Delta$ is invertible as an operator on the range
of $P_{2}$. Equation (\ref{pos}) implies that \[
(P_{2}-D(0)D(z)^{\ast})\Delta^{-2}(P_{2}-D(z)D(0)^{\ast})\leq(\frac{P_{2}-P_{2}F(z)F(z)^{\ast}P_{2}}{1-|z|^{2}}).\]
 Since $D(0)D(z)^{\ast}$ lies in $B(P_{2}(H))$ and has norm strictly
less than $1$ (as $\Vert D(0)\Vert<1$), $P_{2}-D(0)D(z)^{\ast}$
is invertible in $B(P_{2}(H))$ and so, therefore, is $(P_{2}-D(0)D(z)^{\ast})\Delta^{-2}(P_{2}-D(z)D(0)^{\ast})$.
Hence, for each $z\in\mathbb{D}$ there is a $q(z)>0$, such that
$\frac{P_{2}-P_{2}F(z)F(z)^{\ast}P_{2}}{1-|z|^{2}}\geq(P_{2}-D(0)D(z)^{\ast})\Delta^{-2}(P_{2}-D(z)D(0)^{\ast})\geq q(z)P_{2}$.
Thus, \[
P_{2}-P_{2}F(z)F(z)^{\ast}P_{2}\geq(1-|z|^{2})q(z)P_{2},\]
which yields $P_{2}F(z)F(z)^{\ast}P_{2}\leq(1-q(z)(1-|z|^{2}))P_{2}$.
So, if we set $q_{0}(z)=(1-q(z)(1-|z|^{2}))$, we obtain a function
with the desired properties. \end{proof}

We return to our analysis of the special function $F:\mathbb{D}\rightarrow B(E\otimes_{\sigma}H,H)$
defined in equation (\ref{F}).

\begin{lemma} \label{F0}The function $F$ defined by equation (\ref{F})
satisfies:

\begin{enumerate}
\item [(1)] For every $z\in\mathbb{D}$ and $a\in M$, $F(z)(\varphi_{E}(a)\otimes I_{H})=\sigma(a)F(z)$
and $F(z)F(z)^{*}$ commutes with $\sigma(M)$.
\item [(2)] For every $b\in\sigma(M)^{\prime}$, $bF(0)=F(0)(I_{E}\otimes b)$
and $F(0)F(0)^{\ast}\in\mathfrak{Z}(\sigma(M))$. 
\end{enumerate}
\end{lemma}

\begin{proof} Since $F(z)^{\ast}\in E^{\sigma}$ by Lemma~\ref{welldef},
(1) holds. For (2), simply note that $bF(0)(\xi\otimes h)=b\alpha(T_{\xi})(0)h=\alpha(T_{\xi})(0)bh=F(0)(\xi\otimes bh)=F(0)(I_{E}\otimes b)(\xi\otimes h)$,
where we used the fact that for every $X\in H^{\infty}(E)$, $X(0)\in\sigma(M)$.
\end{proof}

\begin{definition} \label{split} Let $\tau$ be the map defined
by equation (\ref{tau}). We say that $\tau(0)$ splits if there are
projections $P_{1},P_{2}$ in $\sigma(M)^{\prime}$ such that

\begin{enumerate}
\item [(i)] $P_{1}+P_{2}=I$,
\item [(ii)] $P_{1}\tau(0)^{*}\tau(0)P_{1}=P_{1}$ and
\item [(iii)] $P_{2}\tau(0)^{*}\tau(0)P_{2}\leq rP_{2}$ for some $r<1$. 
\end{enumerate}
\end{definition}

Note that $\tau(0)=F(0)^{\ast}$ so that, although $F$ depends on
a choice of $\eta\in\mathbb{D}(E^{\sigma})$, $F(0)$ does not. It
follows from Lemma~\ref{F0}, therefore, that $\tau(0)^{\ast}\tau(0)$
lies in the center of $\sigma(M)$, $\mathfrak{Z}(\sigma(M))=\sigma(\mathfrak{Z}(M))$.

Note also that, if the center of $M$, $\mathfrak{Z}(M)$, is an atomic
abelian von Neumann algebra, then $\tau(0)$ always splits. This is
the case, in particular, if $M$ is a factor or if $M=\mathbb{C}^{n}$.
It is also the case, therefore, when $E$ is the correspondence associated
with a (countable) directed graph.

When $\tau(0)$ splits we have the following.

\begin{theorem} \label{strictcontraction} Assume that the left action
map of $M$ on $E$, $\varphi_{E}$, is injective and that $\tau(0)$
splits. Then the map $\tau$ defined in equation (\ref{tau})) maps
$\mathbb{D}(E^{\sigma})$ into itself and satisfies the following
equation \[
(\widehat{\alpha(X)})(\eta^{\ast})=\widehat{X}(\tau(\eta)^{\ast}),\]
 for every $X\in H^{\infty}(E)$ and $\eta\in\mathbb{D}(E^{\sigma})$.
\end{theorem}

\begin{proof} Fix $\eta\in\mathbb{D}(E^{\sigma})$ and let $F$ be
the map defined in (\ref{F}). Since $\tau(0)=F(0)^{\ast}$ splits,
there are projections $P_{1},P_{2}$ as in Definition~\ref{split}.
Using Lemma~\ref{F0}, we see that the conditions of Theorem~\ref{decomp}
are satisfied with $K=E\otimes H$ and $Q_{i}=I_{E}\otimes P_{i}$,
$i=1,2$. Thus, \[
P_{1}F(z)=P_{1}F(z)(I_{E}\otimes P_{1})=P_{1}F(0)(I_{E}\otimes P_{1})=P_{1}F(0)\]
 for all $z\in\mathbb{D}$. Consequently, for all $\xi\in E$, $P_{1}(\sigma\times z\eta^{\ast})(\alpha(T_{\xi}))=P_{1}\sigma(\alpha(T_{\xi})_{0})$
where, for $X\in H^{\infty}(E)$, $X_{0}$ is the image of $X$ under
the conditional expectation onto $\varphi_{\infty}(M)$. Since the
representation $\sigma\times z\eta^{\ast}$ is $w^{\ast}$-continuous
and $\alpha$ is surjective, we have for all $X\in H^{\infty}(E)$,
\[
P_{1}(\sigma\times z\eta^{\ast})(X)=P_{1}\sigma(X_{0}).\]
 In particular, letting $X=T_{\xi}$, we see that $P_{1}(\sigma\times z\eta^{\ast})(T_{\xi})=0$.
Since, for $h\in H$, $(\sigma\times z\eta^{\ast})(T_{\xi})h=P_{1}\eta^{\ast}(\xi\otimes h)=0$
we have $\eta P_{1}=0$. (Recall that $P_{1}\in\sigma(M)^{\prime}$
and, thus, $\eta P_{1}$ is well defined since $E^{\sigma}$ is a
right module over $\sigma(M)^{\prime}$).

Since $\eta$ is arbitrary in $\mathbb{D}(E^{\sigma})$, $E^{\sigma}P_{1}=0$.
If $P_{1}\neq0$, it follows that $E^{\sigma}$ is not full and, using
Lemma~\ref{faithful}, the map $\varphi_{E}$ is not injective, contradicting
our assumption. Thus $P_{1}=0$ and it follows from Theorem~\ref{decomp}
that $\Vert F(z)\Vert<1$ for every $z$. since this holds for all
$\eta\in\mathbb{D}(E^{\sigma})$, the conclusion of the theorem follows.
\end{proof}

Next we show that the map $\tau$ is holomorphic on $\mathbb{D}(E^{\sigma})$.
We view it as a map into $B(H,E\otimes H)$. To be holomorphic is
the same as being Frechet-differentiable. If we use \cite[Theorem 3.17.1]{HP57}
and the fact, proved in Lemma~\ref{welldef}, that $\tau$ is bounded,
it suffices to show that $\tau$ is (G)-differentiable in the sense
of \cite[Definition 3.16.2]{HP57}. But if we apply \cite[Theorem 3.16.1]{HP57},
this means that we have to show that for every $\eta_{0},\eta\in\mathbb{D}(E^{\sigma})$,
the function $G(z):=\tau(\eta_{0}+z\eta)$, defined on $D(\eta,\eta_{0}):=\{ z\in\mathbb{C}||z|<(1-\Vert\eta_{0}\Vert)/\Vert\eta\Vert\}$
is holomorphic in the sense of \cite[Definition 3.10.1]{HP57}.

Since the set of all functionals on $B(H,E\otimes H)$ that are $w^{\ast}$-continuous
is a determining manifold for $B(H,E\otimes H)$ in the sense of \cite[Definition 2.8.2]{HP57},
it suffices to show that for every $w^{\ast}$-continuous functional
$w$, the map $z\mapsto w(\tau(\eta_{0}+z\eta))$ is holomorphic on
$D(\eta,\eta_{0})$. It is enough, in fact, to consider all functionals
of the form $T\mapsto\langle Th,\xi\otimes k\rangle$ for $h,k\in H$
and $\xi$ in the unit ball of $E$.

So we fix $\eta_{0},\eta\in E^{\sigma}$, $h,k\in H$ and $\xi\in E$
with $\Vert\xi\Vert<1$ and write $f(z)=\langle\tau(\eta_{0}+z\eta)h,\xi\otimes k\rangle$
for $z\in D(\eta,\eta_{0})$. We have \[
f(z)=\langle h,\tau(\eta_{0}+z\eta)^{\ast}(\xi\otimes k)\rangle=\langle h,\widehat{\alpha(T_{\xi})}(\eta_{0}^{\ast}+\bar{z}\eta^{\ast})k\rangle.\]
 Note that by Theorem~\ref{evaluationtheorem}, we can write \[
\widehat{\alpha(T_{\xi})}(\eta_{0}^{\ast}+z\eta^{\ast})=A+\sum_{m=1}^{\infty}B((L_{\eta_{0}}^{\ast}+\bar{z}L_{\eta}^{\ast})D)^{m-1}(L_{\eta_{0}}^{\ast}+\bar{z}L_{\eta}^{\ast})C\]
 where $A,B,C,D$ are from some system matrix and the sum converges
in norm. Thus \[
f(z)=\langle A^{\ast}h,k\rangle+\sum_{m=1}^{\infty}\langle C^{\ast}(L_{\eta_{0}}+zL_{\eta})(D^{\ast}(L_{\eta_{0}}+zL_{\eta}))^{m-1}B^{\ast}h,k\rangle\]
 and this function is clearly holomorphic.

We can conclude:

\begin{corollary} \label{holom}The function $\tau$ is a holomorphic
map from $\mathbb{D}(E^{\sigma})$ to its closure. \end{corollary}

\begin{theorem} \label{autom1} Let $E$ be a faithful $W^{\ast}$-correspondence
over $M$, let $\alpha$ be an automorphism of $H^{\infty}(E)$ that
is completely isometric, is a $w^{\ast}$-homeomorphism and leaves
$\varphi_{\infty}(M)$ elementwise fixed, and let $\sigma$ be a faithful
representation of $M$. Write $\tau$ for the transpose of $\alpha$
defined in equation (\ref{tau}) and write $\theta$ for the map associated
similarly with $\alpha^{-1}$. If both $\tau(0)$ and $\theta(0)$
split (as in Definition~\ref{split}) then $\tau$ is a biholomorphic
map of the open unit ball of $E^{\sigma}$, $\tau^{-1}=\theta$, and,
for every $X\in H^{\infty}(E)$, \begin{equation}
(\widehat{\alpha(X)})(\eta^{\ast})=\widehat{X}(\tau(\eta)^{\ast})\;,\;\eta\in\mathbb{D}(E^{\sigma}).\label{implement}\end{equation}

\end{theorem}

\begin{proof} We already know that, under the conditions of the theorem,
both $\tau$ and $\theta$ are holomorphic maps of the open unit ball.
It follows from equation (\ref{tau}) that, for every $\xi\in E$,
$h\in H$ and $\eta\in\mathbb{D}(E^{\sigma})$, $\widehat{\alpha(T_{\xi})}(\eta^{\ast})=\tau(\eta)^{\ast}(\xi\otimes h)$.
But $\tau(\eta)^{\ast}(\xi\otimes h)=\widehat{T_{\xi}}(\tau(\eta)^{\ast})$,
so that equation (\ref{implement}) holds for $T_{\xi}$. It also
holds for $\varphi_{\infty}(a)$, $a\in M$, since $\alpha(\varphi_{\infty}(a))=\varphi_{\infty}(a)$.
Therefore it holds for every $X$ in a $w^{\ast}$-dense subalgebra
of $H^{\infty}(E)$. By the $w^{\ast}$-continuity of $\alpha$, equation
(\ref{implement}) holds for every $X\in H^{\infty}(E)$. Since a
similar claim holds for $\alpha^{-1}$ and $\theta$, we conclude
that for all $X\in H^{\infty}(E)$, $\widehat{X}(\eta^{\ast})=\widehat{\alpha^{-1}(\alpha(X))}(\eta^{\ast})=\widehat{\alpha(X)}(\theta(\eta)^{\ast})=\widehat{X}(\tau(\theta(\eta))^{\ast})$.
Thus $\tau^{-1}=\theta$. \end{proof}

A biholomorphic map $\tau$ is said to \emph{implement} $\alpha$
if equation (\ref{implement}) holds.

\begin{remark} \label{nonid} If $\alpha$ is implemented by $\tau$
in the sense of equation (\ref{implement}), then, writing this equation
when $X=\varphi_{\infty}(a)$, $a\in M$, shows that $\alpha$ leaves
$\varphi_{\infty}(M)$ elementwise fixed. Also, inspecting the proof
of Lemma~\ref{welldef}, one sees that, if $\alpha$ does not have
this property, the map $\tau,$ defined in equation (\ref{tau}) would
map the unit ball of $E^{\sigma}$ into the unit ball of $E^{\pi}$
where $\pi=\sigma\circ\varphi_{\infty}^{-1}\circ\alpha\circ\varphi_{\infty}$.
One can study such automorphisms by studying these maps but the situation
becomes quite complicated, unless one makes a global assumption to
begin with, \emph{vis.}, that $\sigma$ has uniform infinite multiplicity.
In that event, by properties of normal representations of von Neumann
algebras, $\sigma$ and $\pi$ are unitarily equivalent. Say $\pi(\cdot)=u\sigma(\cdot)u^{\ast}$
for some Hilbert space isomorphism from the Hilbert space of $\sigma$
to the Hilbert space of $\pi$. Then it is a straightforward calculation
to see that $E^{\pi}=(I\otimes u)E^{\sigma}u^{\ast}$. It is then
a straightforward matter to incorporate $u$ into our formulas. \end{remark}

As we have remarked before, $\mathbb{D}(E^{\sigma})$ is the unit
ball of a $J^{\ast}$-triple system. It results, therefore, from well-known
theory \cite{lH74} that the biholomorphic maps of $\mathbb{D}(E^{\sigma})$
are determined by Möbius transformations (and {}``isometric multipliers\textquotedblright).
As we shall, however, the Möbius transformations of $\mathbb{D}(E^{\sigma})$
that implement automorphisms of $H^{\infty}(E)$ have to have a special
form: They must be parametrized by {}``central\textquotedblright\ elements
of $\mathbb{D}(E^{\sigma})$ in the sense of the following definition.
(See also Remark 2.1.3 of \cite{BBLS}).

\begin{definition} \label{center}Let $E$ be a $W^{\ast}$-correspondence
over a $W^{\ast}$-algebra $M$. The \emph{center of} $E$\emph{,}
denoted $\mathfrak{Z}(E)$, is the set of $\xi\in E$ such that $a\xi=\xi a$
for all $a\in M$. \end{definition}

\begin{lemma} \label{propcenter}

\begin{enumerate}
\item [(1)] The center $\mathfrak{Z}(E)$ of a $W^{\ast}$-correspondence
$E$ over $M$ is a $W^{\ast}$-correspondence over the center $\mathfrak{Z}(M)$
of $M$.
\item [(2)] Let $\sigma$ be a faithful normal representation of $M$
on the Hilbert space $\mathcal{E}$, and for $\xi\in E$, define $\Phi(\xi):=L_{\xi}$
where $L_{\xi}$ maps $\mathcal{E}$ to $E\otimes\mathcal{E}$ via
the formula $L_{\xi}(h)=\xi\otimes h$. Then the pair $(\sigma,\Phi)$
defines an isomorphism of $\mathfrak{Z}(E)$ onto $\mathfrak{Z}(E^{\sigma})$
in the sense of Definition~\ref{isomorph}. (Here, $\mathfrak{Z}(E)$
is a correspondence over $\mathfrak{Z}(M)$ and $\mathfrak{Z}(E^{\sigma})$
is a correspondence over $\mathfrak{Z}(\sigma(M)^{\prime})=\mathfrak{Z}(\sigma(M))=\sigma(\mathfrak{Z}(M))$).
\item [(3)] Given a faithful representation $\sigma$ of $M$ on the Hilbert
space $\mathcal{E}$ and $\gamma\in\mathbb{D}(E^{\sigma})$, then
$\gamma$ lies in the center of $E^{\sigma}$ if and only if the representation
$\sigma\times\gamma^{\ast}$ maps $H^{\infty}(E)$ into $\sigma(M)$. 
\end{enumerate}
\end{lemma}

\begin{proof} It is clear that $\mathfrak{Z}(E)$ is a bimodule over
$\mathfrak{Z}(M)$ and, to prove (1), we need only show that the inner
product of two elements in $\mathfrak{Z}(E)$ lies in $\mathfrak{Z}(M)$.
For $a\in M$, $\xi_{1},\xi_{2}\in\mathfrak{Z}(E)$ we have \[
a\langle\xi_{1},\xi_{2}\rangle=\langle\xi_{1}a^{\ast},\xi_{2}\rangle=\langle a^{\ast}\xi_{1},\xi_{2}\rangle=\langle\xi_{1},a\xi_{2}\rangle=\langle\xi_{1},\xi_{2}a\rangle=\langle\xi_{1},\xi_{2}\rangle a.\]
 Hence the inner product lies in the center of $M$, proving (1).
We fix a faithful representation $\sigma$ of $M$ on $\mathcal{E}$.
For $\xi\in\mathfrak{Z}(E)$, $a\in M$ and $h\in\mathcal{E}$ we
have $L_{\xi}\sigma(a)h=\xi\otimes_{\sigma}\sigma(a)h=\xi a\otimes h=a\xi\otimes h=(a\otimes I)L_{\xi}h$.
Hence, $L_{\xi}\in E^{\sigma}$. Given $b\in\sigma(M)^{\prime}$ and
$h\in\mathcal{E}$ we have $L_{\xi}bh=\xi\otimes bh=(I_{E}\otimes b)L_{\xi}h$.
Thus $L_{\xi}$ lies in $\mathfrak{Z}(E^{\sigma})$.

For $\xi\in\mathfrak{Z}(E)$, $a,b\in\mathfrak{Z}(M)$, and $h\in\mathcal{E}$,
$L_{a\xi b}h=a\xi b\otimes h=\xi ab\otimes h=\xi\otimes\sigma(a)\sigma(b)h=(I\otimes\sigma(a))L_{\xi}\sigma(b)h$
hence \[
\Phi(a\xi b)=\sigma(a)\Phi(\xi)\sigma(b).\]
 For $\xi_{1},\xi_{2}\in\mathfrak{Z}(E)$ we have $L_{\xi_{1}}^{\ast}L_{\xi_{2}}=\sigma(\langle\xi_{1},\xi_{2}\rangle)$.
Therefore the pair $(\sigma,\Phi)$ is an isomorphism of $\mathfrak{Z}(E)$
into $\mathfrak{Z}(E^{\sigma})$.

To prove that the map $\Phi$ is onto, fix an $\eta\in\mathfrak{Z}(E^{\sigma})$.
Then, $\eta$ is a map from $\mathcal{E}$ to $E\otimes_{\sigma}\mathcal{E}$
satisfying \begin{equation}
\eta\sigma(a)=(a\otimes I)\eta\label{a}\end{equation}
 and \begin{equation}
\eta b=(I\otimes b)\eta,\label{b}\end{equation}
 for $a\in M$ and $b\in\sigma(M)^{\prime}$. Define the map $\psi:E\rightarrow B(\mathcal{E})$
by $\psi(\zeta)=\eta^{\ast}L_{\zeta}$ and note that for $b\in\sigma(M)^{\prime}$
and $h\in\mathcal{E}$, $\eta^{\ast}L_{\zeta}bh=\eta^{\ast}(\zeta\otimes bh)=\eta^{\ast}(I\otimes b)L_{\zeta}h$.
Using (\ref{b}) the latter is equal to $b\eta^{\ast}L_{\zeta}h$.
Hence $\psi(\zeta)$ lies in $\sigma(M)$. Also $\psi(\zeta a)=\psi(\zeta)\sigma(a)$
for all $a\in M$ and it then follows from the self duality of $E$
that there is an $\xi\in E$ with $\langle\xi,\zeta\rangle=\sigma^{-1}(\psi(\zeta))$.
Thus, for all $\zeta\in E$, $L_{\xi}^{\ast}L_{\zeta}=\sigma(\langle\xi,\zeta\rangle)=\eta^{\ast}L_{\zeta}$
and we conclude that $\eta=L_{\xi}$.

It follows from (\ref{a}) that, for all $a\in M$, $L_{\xi a}=\eta\sigma(a)=(a\otimes I)\eta=L_{a\xi}$,
showing that $\xi$ lies in $\mathfrak{Z}(E)$.

Finally, to prove (3), fix an $\eta\in\mathbb{D}(E^{\sigma})$ and
write $(T,\sigma)$ for the covariant pair associated with $\sigma\times\eta^{\ast}$
(so that, $\tilde{T}=\eta^{\ast}$). Then the representation maps
$H^{\infty}(E)$ into $\sigma(M)$ if and only if, for each $\xi\in E$,
$T(\xi)\in\sigma(M)$. This holds iff, for all $b\in\sigma(M)^{\prime}$,
$\xi\in E$ and $h\in\mathcal{E}$, $\tilde{T}(I_{\mathcal{E}}\otimes b)(\xi\otimes h)=T(\xi)bh=bT(\xi)h=b\tilde{T}(\xi\otimes h)$;
that is, if and only if $\tilde{T}(I_{\mathcal{E}}\otimes b)=b\tilde{T}$
for every $b\in\sigma(M)^{\prime}$. But the last statement says that
$\eta$ lies in the center of $E^{\sigma}$. \end{proof}

The following example may help to show that the center of a correspondence
is much less {}``inert\textquotedblright\ than the center of a
von Neumann algebra.

\begin{example} \label{IntertwinerHmodule}Let $M$ be a von Neumann
algebra and let $\alpha$ be an endomorphism of $M$. Then we obtain
a $W^{\ast}$-correspondence over $M$, denoted $_{\alpha}M$, by
taking $M$ with its usual right action and inner product give by
the formula, $\langle\xi,\eta\rangle=\xi^{\ast}\eta$ and by letting
$\alpha$ implement the left action. Then an element $\xi$ in $_{\alpha}M$
lies in the center of $_{\alpha}M$ if and only if $\xi$ intertwines
$\alpha$ and the identity endomorphism; i.e., $\xi\in\mathfrak{Z}(_{\alpha}M)$
if and only if $\alpha(a)\xi=\xi a$ for all $a\in M$. $\mathfrak{Z}(_{\alpha}M)$
is a much studied object in the literature and the preceding lemma
spells out some of its important elementary properties. \end{example}

Our goal now is to develop the properties of Möbius transformations
of $\mathbb{D}(E^{\sigma})$ and to identify those that implement
automorphisms of $H^{\infty}(E)$. To this end, fix a faithful representation
$\sigma$ of $M$ on a Hilbert space $\mathcal{E}$. Set $N=\sigma(M)^{\prime}$,
write $K=\mathcal{E}\oplus(E\otimes_{\sigma}\mathcal{E})$, and define
the (necessarily faithful) representation $\rho$ of $N$ on $K$
by the formula \[
\rho(S)=\left(\begin{array}{cc}
S & 0\\
0 & I\otimes S\end{array}\right),\;\; S\in N.\]
 For $\gamma\in\mathbb{D}(E^{\sigma})$ we set $\Delta_{\gamma}:=(I_{\mathcal{E}}-\gamma^{\ast}\gamma)^{1/2}$
- an element in $B(\mathcal{E})$ - and $\Delta_{\gamma\ast}:=(I_{E\otimes\mathcal{E}}-\gamma\gamma^{\ast})^{1/2}$
- an element in $B(E\otimes\mathcal{E})$. When $\gamma$ is understood,
then we shall simply write $\Delta$ for $\Delta_{\gamma}$ and $\Delta_{\ast}$
for $\Delta_{\gamma^{\ast}}$. Given $\gamma\in\mathbb{D}(E^{\sigma})$
we define the map $g_{\gamma}$ on $\mathbb{D}(E^{\sigma})^{\ast}$
by the formula, \begin{equation}
g_{\gamma}(z^{\ast})=\Delta_{\gamma}(I-z^{\ast}\gamma)^{-1}(\gamma^{\ast}-z^{\ast})\Delta_{\gamma^{\ast}}^{-1},\label{ggamma}\end{equation}
 $z\in\mathbb{D}(E^{\sigma})$. Then $g_{\gamma}$ is a biholomorphic
automorphism of $\mathbb{D}(E^{\sigma})^{\ast}$ that maps $0$ to
$\gamma^{\ast}$ and $\gamma^{\ast}$ to $0$. Further, $g_{\gamma}^{2}=id$,
and every biholomorphic map $g$ of $\mathbb{D}(E^{\sigma})^{\ast}$
is of the form \[
g=w\circ g_{\gamma}\]
 where $w$ is an isometry on $(E^{\sigma})^{\ast}$ and $\gamma^{\ast}=w^{-1}g(0)$
\cite{lH74}. When $\gamma$ lies in the center of $E^{\sigma}$,
we see that $g_{\gamma}$ maps the center onto itself and it follows
that every biholomorphic automorphism of the open unit ball of $(E^{\sigma})^{\ast}$
that preserves the center is of the form \[
g=w\circ g_{\gamma}\]
 where $\gamma$ lies in the center and $w$ is an isometry on $(E^{\sigma})^{\ast}$
that preserves the center.

If $z\in\mathbb{D}(E^{\sigma})$, then the series $\sum_{n=0}^{\infty}(z^{\ast}\gamma)^{n}$
converges in norm to the operator in $N$, $(I-z^{\ast}\gamma)^{-1}=\sum_{n=0}^{\infty}(z^{\ast}\gamma)^{n}$.
One easily calculates, then, that \[
g_{\gamma}(z^{\ast})=\Delta\gamma^{\ast}\Delta_{\ast}^{-1}-\Delta(I-z^{\ast}\gamma)^{-1}z^{\ast}\Delta_{\ast}.\]
 Recall that the equation $U(z\otimes h)=z(h)$ defines a Hilbert
space isomorphism $U:E^{\sigma}\otimes\mathcal{E}\rightarrow E\otimes\mathcal{E}$
\cite[p. 369]{MS03}. Consequently, as maps on $\mathcal{E}$, $UL_{z}=z$
and $z^{\ast}=L_{z}^{\ast}U^{\ast}$. Thus we may write \[
g_{\gamma}(z^{\ast})=\Delta\gamma^{\ast}\Delta_{\ast}^{-1}-\Delta(I-L_{z}^{\ast}U^{\ast}\gamma)^{-1}L_{z}^{\ast}U^{\ast}\Delta_{\ast}.\]

We write $K_{1}=E\otimes_{\sigma}\mathcal{E}$ for the second summand
in $K=\mathcal{E}\oplus(E\otimes_{\sigma}\mathcal{E})$ and we let
$q_{1}$ denote the projection from $K$ onto $K_{1}$. Likewise,
we set $K_{2}=\mathcal{E}$ with projection $q_{2}$. Corresponding
to the direct sum decomposition, we define $V$ by the formula \begin{equation}
V:=\left(\begin{array}{cc}
\Delta\gamma^{\ast}\Delta_{\ast}^{-1} & -\Delta\\
U^{\ast}\Delta_{\ast} & U^{\ast}\gamma\end{array}\right):\left(\begin{array}{c}
K_{1}\\
\mathcal{E}\end{array}\right)\rightarrow\left(\begin{array}{c}
K_{2}\\
E^{\sigma}\otimes\mathcal{E}\end{array}\right).\label{V}\end{equation}
 If we calculate $VV^{\ast}$, we find that the off diagonal terms
vanish and the terms on the diagonal are $\Delta\gamma^{\ast}\Delta_{\ast}^{-2}\gamma\Delta+\Delta^{2}$
and $U^{\ast}(\Delta_{\ast}^{2}+\gamma\gamma^{\ast})U$. Since $\Delta_{\ast}^{2}+\gamma\gamma^{\ast}=I_{E\otimes\mathcal{E}}$,
the latter expression is $U^{\ast}U=I_{E^{\sigma}\otimes\mathcal{E}}=q_{2}$.
For the first expression, we note that $\gamma^{\ast}\Delta_{\ast}^{-2}\gamma=\gamma^{\ast}(I-\gamma\gamma^{\ast})^{-1}\gamma=(I-\gamma^{\ast}\gamma)^{-1}-1$
and $\Delta\gamma^{\ast}\Delta_{\ast}^{-2}\gamma\Delta+\Delta^{2}=\Delta((I-\gamma^{\ast}\gamma)^{-1}-I)\Delta+\Delta^{2}=I_{\mathcal{E}}$.
This shows that $V$ is a coisometry. Similar computations show that
it is, in fact, a unitary operator. Thus $V$ is a transfer operator.

We want to apply Theorem~\ref{evaluationtheorem} to obtain an element
$X\in H^{\infty}(E)$ with $\widehat{X}(\eta^{\ast})=g_{\gamma}(\eta^{\ast})$,
for $\eta\in\mathbb{D}(E^{\sigma})$. To do this, we first let $F$
be the correspondence $E^{\sigma}$ and then $F^{\rho}$ is a correspondence
over $\rho(N)^{\prime}$. In order to apply Theorem~\ref{evaluationtheorem}
we let $M$, in that theorem, be the von Neumann algebra $\rho(N)^{\prime}$
and let $\sigma$ there be the identity representation of $\rho(N)^{\prime}$
on $K$ (so that $\mathcal{E}$ there is $K$). $E$ in that theorem
will be $F^{\rho}$ and $N$ there (the commutant of $\sigma(M)$)
will be $\rho(N)$. The representation $\tau$ of $N$ then will be
the map $\rho^{-1}$ of $\rho(N)$ on $\mathcal{E}$ (so that $\mathcal{E}$
will play the role of $H$ there). Also, $q_{1}$ will be as above.
We set $A=\Delta\gamma^{\ast}\Delta_{\ast}^{-1}$, $B=-\Delta$, $C=U^{\ast}\Delta_{\ast}$
and $D=U^{\ast}\gamma$. These $A,B,C$ and $D$ give rise to the
matricial operator $V$ of equation (\ref{V}). In order to show that
the assumptions of Theorem~\ref{evaluationtheorem} are satisfied,
we have to show that these operators ($A,B,C$ and $D$) all have
the required intertwining properties. (Note that we have already checked
that $V$ is a unitary operator).

The required intertwining properties are:

\begin{itemize}
\item [(a)] $A=\Delta\gamma^{*}\Delta_{*}^{-1}$ lies in $q_{2}\rho(N)^{\prime}q_{1}$.
\item [(b)] $B=-\Delta$ lies in $N^{\prime}$.
\item [(c)] For every $S\in N$, $U^{*}\Delta_{*}(I_{E}\otimes S)=(S\otimes I_{\mathcal{E}})U^{*}\Delta_{*}$
on $E\otimes\mathcal{E}$.
\item [(d)] For every $S\in N$, $U^{*}\gamma S=(I_{E}\otimes S)U^{*}\gamma$
on $\mathcal{E}$. 
\end{itemize}
Indeed, recall that $\gamma$ lies in the center of $E^{\sigma}$
and, thus, for $S\in N$, $\gamma S=(I\otimes S)\gamma$. Therefore
$\Delta$ commutes with $N$ and $\Delta_{*}$ commutes with $I\otimes S$
for $S\in N$. This implies (a) and (b). Recall that, for $h\in\mathcal{E}$,
$U^{*}\gamma h=\gamma\otimes h$ and, thus, $U^{*}\gamma Sh=\gamma\otimes Sh=(I\otimes S)(\gamma\otimes h)=(I\otimes S)U^{*}\gamma h$
proving (d). For (c), it suffices to note that $U(S\otimes I)U^{*}=I\otimes S$
and $\Delta_{*}$ commutes with $I\otimes S$ for all $S\in N$.

We can now apply Theorem~\ref{evaluationtheorem}. Since $F^{\rho}$
plays the role of $E$ in that theorem and the identity representation
of $\rho(N)^{\prime}$, $id$, plays the role of $\sigma$, $E^{\sigma}$
in that theorem is replaced by $(F^{\rho})^{id}$ which, by the duality
theorem \cite[Theorem 3.6 ]{MS03} is isomorphic to $F=E^{\sigma}$.
We therefore conclude:

\begin{lemma} \label{g}For every $\gamma\in\mathbb{D}(\mathfrak{Z}(E^{\sigma}))$,
there is an $X$ in $H^{\infty}(F^{\rho})$ with $\Vert X\Vert\leq1$
such that, for all $z\in\mathbb{D}(E^{\sigma})$, $\widehat{X}(z^{\ast})=g_{\gamma}(z^{\ast}).$
\end{lemma}

Note that $g_{\gamma}(z^{\ast})$ is an operator from $E\otimes\mathcal{E}$
into $\mathcal{E}$ and can be viewed as an operator in $B(K)$ which
is where the values of $X$, as an element of $H^{\infty}(F^{\rho})$,
lie.

We can now use \cite[Theorem 5.3]{MS03} to prove the following.

\begin{corollary} \label{kerg} Fix $\gamma\in\mathbb{D}(\mathfrak{Z}(E^{\sigma}))$
as above. Then, for every $z_{1},z_{2},\ldots,z_{k}$ in $\mathbb{D}(E^{\sigma})$,
the map on $M_{k}(\sigma(M)^{\prime})$ defined by the $k\times k$
matrix \[
((id-\theta_{g_{\gamma}(z_{i}^{\ast})^{\ast},g_{\gamma}(z_{j}^{\ast})^{\ast}})\circ(id-\theta_{z_{i},z_{j}})^{-1})\]
 is completely positive. \end{corollary}

\begin{proof} Applying \cite[Theorem 5.3]{MS03} to $X$ of Lemma~\ref{g},
we get the complete positivity of the map defined by the matrix \[
((I-Ad(g_{\gamma}(z_{i}^{\ast}),g_{\gamma}(z_{j}^{\ast})))\circ(id-\theta_{z_{i},z_{j}})^{-1}).\]
 But note that, for every $b\in\sigma(M)^{\prime}$, $Ad(g_{\gamma}(z_{i}^{\ast}),g_{\gamma}(z_{j}^{\ast}))(\rho(b))=g_{\gamma}(z_{i}^{\ast})\rho(b)g_{\gamma}(z_{j}^{\ast})^{\ast}=\langle g_{\gamma}(z_{i}^{\ast})^{\ast},bg_{\gamma}(z_{j}^{\ast})^{\ast}\rangle=\theta_{g_{\gamma}(z_{i}^{\ast})^{\ast},g_{\gamma}(z_{j}^{\ast})^{\ast}}(b)$.
\end{proof}

\begin{corollary} \label{comp} Let $Z:\mathbb{D}(E^{\sigma})^{\ast}\rightarrow B(\mathcal{E})$
be a Schur class operator function and let $\gamma$ be in $\mathbb{D}(\mathfrak{Z}(E^{\sigma}))$.
Then the function $Z_{\gamma}:\mathbb{D}((E^{\sigma})^{\ast})\rightarrow B(\mathcal{E})$
defined by \[
Z_{\gamma}(\eta^{\ast})=Z(g_{\gamma}(\eta^{\ast}))\]
 is also a Schur class operator function. \end{corollary}

\begin{proof} For every $\eta_{i},\eta_{j}$ in $\mathbb{D}(E^{\sigma})$
we have $(id-Ad(Z(g_{\gamma}(\eta_{i}^{*})),Z(g_{\gamma}(\eta_{j}^{*}))))\circ(id-\theta_{\eta_{i},\eta_{j}})^{-1}=((id-Ad(Z(g_{\gamma}(\eta_{i}^{*})),Z(g_{\gamma}(\eta_{j}^{*}))))\circ(id-\theta_{g_{\gamma}(\eta_{i}^{*})^{*},g_{\gamma}(\eta_{j}^{*})^{*}})^{-1})\circ(id-\theta_{g_{\gamma}(\eta_{i}^{*})^{*},g_{\gamma}(\eta_{j}^{*})^{*}})\circ(id-\theta_{\eta_{i},\eta_{j}})^{-1})$.
Hence the map associated with $Z_{\gamma}$ is a composition of two
completely positive maps and is, therefore, completely positive. \end{proof}

For the statement of the next lemma, recall from \cite[end of Section 2]{MS03}
that every $X\in H^{\infty}(E)$ has a {}``Fourier series\char`\"{}
expansion given by a sequence of {}``Fourier coefficient operators\char`\"{}
$\{\mathbb{E}_{j}\}$. (In \cite{MS03} we wrote $\{\Phi_{j}\}$ for
this sequence). Each map $\mathbb{E}_{j}$ is completely contractive,
$w^{\ast}$-continuous and $\mathbb{E}_{j}(T_{\xi_{1}}T_{\xi_{2}}\cdots T_{\xi_{k}})=T_{\xi_{1}}T_{\xi_{2}}\cdots T_{\xi_{k}}$
if $j=k$ and is zero otherwise. The Cesaro means of the {}``Fourier
series\char`\"{} of $X$ converge to $X$ in the $w^{\ast}$-topology.

\begin{lemma} \label{kersigma} Let $\sigma$ be a normal, faithful,
representation of $M$ on a Hilbert space $H$ and let $K(\sigma)$
denote the kernel of the map $X\rightarrow\widehat{X}$ defined in
equation (\ref{KernalSigma}).

\begin{enumerate}
\item [(i)] $K(\sigma)\subseteq\{ X\in H^{\infty}(E)\;|\;\mathbb{E}_{0}(X)=\mathbb{E}_{1}(X)=0\}$.
\item [(ii)] If, for every $k\in\mathbb{N}$, $\vee\{(\eta^{\otimes k})(H)\;|\;\eta\in\mathbb{D}(E^{\sigma})\}=E^{\otimes k}\otimes H$,
then $K(\sigma)=\{0\}$.
\item [(iii)] Every completely isometric automorphism $\alpha$ of $H^{\infty}(E)$
that is a $w^{*}$-homeomorphism and is implemented by a biholomorphic
map of $\mathbb{D}(E^{\sigma})$ in the sense of (\ref{implement})
leaves $K(\sigma)$ invariant. In particular, $K(\sigma)$ is invariant
under the action of the gauge group and, thus, under the maps $\mathbb{E}_{k}$,
$k\geq0$. 
\end{enumerate}
\end{lemma}

\begin{proof} Write $C_{1}$ for $\{ X\in H^{\infty}(E)\;|\;\mathbb{E}_{0}(X)=\mathbb{E}_{1}(X)=0\}$.
Then for every $X\in H^{\infty}(E)$, $X=\mathbb{E}_{0}(X)+\mathbb{E}_{1}(X)+X_{1}$
where $X_{1}\in C_{1}$. Note that for every $\eta\in\mathbb{D}(E^{\sigma})$,
every $0<t\leq1$ and every $k\geq0$, $\mathbb{E}_{k}(X)((t\eta)^{*})=t^{k}\mathbb{E}(X)(\eta^{*})$.
Thus, for $X\in K(\sigma)$, $0=X((t\eta)^{*})=\mathbb{E}_{0}(X)(\eta^{*})+t\mathbb{E}_{1}(X)(\eta^{*})+t^{2}S$
where $S$ is some bounded operator on $H$. Since this holds for
every $0<t\leq1$, we have (by differentiation) $\mathbb{E}_{0}(X)=0$
and $\mathbb{E}_{1}(X)(\eta^{*})=0$ for all $\eta\in\mathbb{D}(E^{\sigma})$.
Write $\mathbb{E}_{1}(X)=T_{\xi}$ (for some $\xi\in E$). Then, for
all $h\in H$ and $\eta\in\mathbb{D}(E^{\sigma})$, $0=\mathbb{E}_{1}(X)(\eta^{*})h=\eta^{*}(\xi\otimes h)$.
Since $\vee\{\eta(H)|\;\eta\in\mathbb{D}(E^{\sigma})\}=E\otimes H$
(\cite[Lemma 3.5]{MS03}), we find that $\xi\otimes h=0$ for all
$h\in H$. Since $E$ is faithful, this implies that $\xi=0$, completing
the proof of (i).

We can also write $0=X((t\eta)^{*})=\mathbb{E}_{0}(X)(\eta^{*})+t\mathbb{E}_{1}(X)(\eta^{*})+\cdots+t^{k}\mathbb{E}_{k}(X)(\eta^{*})+t^{k+1}S$
and conclude that $\mathbb{E}_{j}(X)(\eta^{*})=0$ for all $j\leq k$.
We can then continue as above but to be able to conclude that $\mathbb{E}_{k}(X)=0$
we need the condition in part (ii) (to replace the use of \cite[Lemma 3.5]{MS03}
in the argument above).

To prove (iii), note that the invariance of $K(\sigma)$ under an
automorphism $\alpha$ as in (iii) follows from (\ref{implement}).
The invariance under the gauge group (and under $\mathbb{E}_{k}$)
is then immediate. \end{proof}

The following proposition is obvious if $K(\sigma)=\{0\}$. But, in
fact, it holds for every faithful, normal representation $\sigma$.
The argument uses an idea from \cite[Proof of Theorem 4.11]{DP98}.

\begin{proposition} \label{uniqueness} Let $\sigma$ be a faithful,
normal representation of $M$ and let $\alpha,\beta$ be two homomorphisms
of $H^{\infty}(E)$ into itself such that $\beta$ is completely isometric,
surjective and a $w^{\ast}$-homeomorphism, while $\alpha$ is completely
contractive and $w^{\ast}$-continuous. Suppose they satisfy the equation
\[
\widehat{\alpha(X)}(\eta^{\ast})=\widehat{\beta(X)}(\eta^{\ast})\]
 for all $X\in H^{\infty}(E)$ and $\eta\in\mathbb{D}(E^{\sigma})$.
Then $\alpha=\beta$. \end{proposition}

\begin{proof} It is clearly enough to assume $\beta=id$ and $\widehat{\alpha(X)}(\eta^{*})=\widehat{X}(\eta^{*})$.
Note that $\alpha$, viewed as a representation of $H^{\infty}(E)$
on $\mathcal{F}(E)\otimes_{\sigma}H$ (whose restriction to $\varphi_{\infty}(M)$
is $\varphi_{\infty}(\cdot)\otimes I_{H}$), can be written as $(\varphi_{\infty}(\cdot)\otimes I_{H})\times\zeta^{*}$
for some $\zeta$ in the closed unit ball of the $\varphi_{\infty}(\cdot)\otimes I_{H}$-dual
of $E$. Thus, for $k\in\mathcal{F}(E)\otimes_{\sigma}H$, $\alpha(T_{\xi})k=(\zeta^{*})(\xi\otimes k)$
and $\Vert\alpha(T_{\xi})k\Vert\leq\Vert\xi\otimes k\Vert=\Vert T_{\xi}k\Vert$.

Fix $h\in H\ $viewed as the zero$^{th}$ summand of $\mathcal{F}(E)\otimes_{\sigma}H$.
Then for every $\xi\in E$, \[
\Vert\alpha(T_{\xi})h\Vert\leq\Vert T_{\xi}h\Vert.\]
 By construction $\alpha(T_{\xi})-T_{\xi}\in K(\sigma)$. But also,
by Lemma~\ref{kersigma}(i), for every $X\in K(\sigma)$, $Xh$ is
orthogonal to $T_{\xi}h$. Thus \[
\Vert\alpha(T_{\xi})h\Vert^{2}=\Vert(\alpha(T_{\xi})-T_{\xi})h\Vert^{2}+\Vert T_{\xi}h\Vert^{2}\geq\Vert T_{\xi}h\Vert^{2}.\]
 We conclude that for every $h\in H$, $(\alpha(T_{\xi})-T_{\xi})h=0$.
It follows that $\alpha(T_{\xi})=T_{\xi}$ for all $\xi\in E$. Since
$\alpha$ is a $w^{\ast}$-continuous homomorphism, $\alpha(X)=X$
for all $X\in H^{\infty}(E)$. \end{proof}

The following lemma will prove very useful when we deal with a representation
$\sigma$ for which $K(\sigma)\neq\{0\}$. It relates the $\sigma$-dual
with the $\pi$-dual where $\pi$ is the representation defined in
the proof of Lemma~\ref{1to1} (for which $K(\pi)=\{0\}$). 

\begin{lemma} \label{lifting} Let $\sigma$ be a faithful representation
of $M$ on $H$ and $\pi$ be the representation $\varphi_{\infty}\otimes I_{H}$
of $M$ on $K:=\mathcal{F}(E)\otimes H$. Let $\psi:\sigma(M)^{\prime}\rightarrow(\varphi_{\infty}(M)\otimes I_{H})^{\prime}$
be defined by $\psi(b)=I_{E}\otimes b$ and let $\Psi:E^{\sigma}\rightarrow E^{\pi}$
be defined by $\Psi(\eta)=I_{\mathcal{F}(E)}\otimes\eta$. Then we
have the following.

\begin{enumerate}
\item [(1)] The pair $(\psi,\Psi)$ is an isomorphism of $E^{\sigma}$
into (not necessarily onto) $E^{\pi}$ satisfying \[
\Psi(\eta)P_{H}=P_{E\otimes H}\Psi(\eta)=\eta\;,\;\eta\in E^{\sigma}\]
 where $P_{H}$ is the projection from $K$ to $H$ (viewed as a subspace)
and $P_{E\otimes H}$ is the projection of $E\otimes K$ onto $E\otimes H$.

\item [(2)] For every $X\in H^{\infty}(E)$ and $\zeta\in E^{\pi}$ that
satisfies $\zeta P_{H}=P_{E\otimes H}\zeta$, we have $\zeta|H\in E^{\sigma}$
and the restriction of $\widehat{X}(\zeta^{*})$ to $H$ (viewed as
a summand of $\mathcal{F}(E)\otimes H=H\oplus E\otimes H\oplus\cdots$)
is $\widehat{X}((\zeta|H)^{*})$.
\item [(3)] There is an isomorphism $\Phi$ of $\mathfrak{Z}(E^{\sigma})$
onto $\mathfrak{Z}(E^{\pi})$ satisfying \[
\Phi(\gamma)P_{H}=P_{E\otimes H}\Phi(\gamma)=\gamma\;,\;\gamma\in\mathfrak{Z}(E^{\sigma}).\]

\item [(4)] For $\eta\in E^{\sigma}$ and $\gamma\in\mathfrak{Z}(E^{\sigma})$,
\[
g_{\Phi(\gamma)}(\Psi(\eta)^{*})P_{E\otimes H}=P_{H}g_{\Phi(\gamma)}(\Psi(\eta)^{*})=g_{\gamma}(\eta^{*}).\]

\end{enumerate}
\end{lemma}

\begin{proof} It is clear that $\psi$ is indeed an isomorphism into
$(\varphi_{\infty}(M)\otimes I_{H})^{\prime}$. Note that it follows
from the intertwining property of $\eta\in E^{\sigma}$ that $\Psi(\eta)$
is a well defined bounded operator. To show that $\Psi$ maps $E^{\sigma}$
into $E^{\pi}$, fix $\eta\in E^{\sigma}$, $\theta\otimes h\in\mathcal{F}(E)\otimes H$
and $a\in M$ and compute $(I_{\mathcal{F}(E)}\otimes\eta)\pi(a)(\theta\otimes h)=(I_{\mathcal{F}(E)}\otimes\eta)(\varphi_{\infty}(a)\theta\otimes h)=\varphi_{\infty}(a)\theta\otimes\eta(h),$
where we view $\mathcal{F}(E)\otimes E$ as the subspace of $\mathcal{F}(E)$
consisting of all the positive tensor powers of $E$. But the last
expression is equal to $(\varphi_{\infty}(a)\otimes I_{H})(I_{\mathcal{F}(E)}\otimes\eta)(\theta\otimes h)$,
showing that $\Psi(\eta)\in E^{\pi}$.

To show that the map is a bimodule map, fix $\eta\in E^{\sigma}$,
$b,c\in\sigma(M)^{\prime}$ and $\theta\otimes h\in\mathcal{F}(E)\otimes H$.
Then $\Psi(c\eta b)(\theta\otimes h)=\theta\otimes(c\eta b)h=\theta\otimes(I_{E}\otimes c)\eta bh=\psi(c)(\theta\otimes\eta bh)=\psi(c)\Psi(\eta)(\theta\otimes bh)=\psi(c)\Psi(\eta)\psi(b)(\theta\otimes h)$,
proving that the image of $\Psi$ lies in $E^{\pi}$. Regarding the
inner product, we have: $\langle\Psi(\eta_{1}),\Psi(\eta_{2})\rangle=\Psi(\eta_{1})^{\ast}\Psi(\eta_{2})=(I_{\mathcal{F}(E)}\otimes\eta_{1})^{\ast}(I_{\mathcal{F}(E)}\otimes\eta_{2})=(I_{\mathcal{F}(E)}\otimes\eta_{1}^{\ast}\eta_{2})=\psi(\langle\eta_{1},\eta_{2}\rangle)$
for all $\eta_{1},\eta_{2}\in E^{\sigma}$. Thus $(\psi,\Psi)$ is
an isomorphism of $E^{\sigma}$ into $E^{\pi}$. The proof of the
equation $\Psi(\eta)P_{H}=P_{E\otimes H}\Psi(\eta)=\eta$ for $\eta\in E^{\sigma}$
is easy. This proves (1). 

To prove (2), let $\zeta\in E^{\pi}$ satisfy $\zeta P_{H}=P_{E\otimes H}\zeta$
and fix $a\in M$ and $h\in H$. Then $(\zeta|H)\sigma(a)h=\zeta(\varphi_{\infty}(a)\otimes I_{H})h=(\varphi_{E}(a)\otimes I_{K})P_{E\otimes H}\zeta h=(\varphi_{E}(a)\otimes I_{H})(\zeta|H)h$.
Thus, $\zeta|H\in E^{\sigma}$. To prove that $\widehat{X}((\zeta|H)^{*})=\widehat{X}(\zeta^{*})|H$,
let, first, consider $X=\varphi_{\infty}(a)$ for $a\in M$. Then
$\widehat{X}(\zeta^{*})=\varphi_{\infty}(a)\otimes I_{H}$ and $\widehat{X}((\eta|H)^{*})=\sigma(a)$
and (2) holds in this case. Take $X=T_{\xi}$ for some $\xi\in E$.
Then, for $h\in H\subseteq\mathcal{F}(E)\otimes H$, $\widehat{X}(\zeta^{*})h=\zeta^{*}(\xi\otimes h)=(\zeta|H)^{*}(\xi\otimes h)=\widehat{X}((\zeta|H)^{*})h$.
In particular, we see that $H$ is invariant for all $\widehat{X}(\zeta^{*})$
where $X$ runs over a set of generators. Thus, $H$ is invariant
under $\widehat{X}(\zeta^{*})$ for all $X\in H^{\infty}(E)$ and
(2) holds for all $X$'s in a $w^{*}$-dense subalgebra of $H^{\infty}(E)$.
Since the map $X\mapsto\widehat{X}(\zeta^{*})$ is $w^{*}$-continuous,
we are done.

To prove (3), recall from Lemma~\ref{propcenter} (2) that both $\mathfrak{Z}(E^{\sigma})$
and $\mathfrak{Z}(E^{\pi})$ are isomorphic to $\mathfrak{Z}(E)$.
Combining these two isomorphisms, we get $\Phi$. More precisely,
every $\eta\in\mathfrak{Z}(E^{\sigma})$ is equal to $L_{\xi}$ for
some $\xi\in\mathfrak{Z}(E)$ (that is, $\eta(h)=\xi\otimes h$, $h\in H$).
Then we set $\Phi(\eta)k=\xi\otimes k$ for $k\in K=\mathcal{F}(E)\otimes H$.
The equation $\Phi(\gamma)P_{H}=P_{E\otimes H}\Phi(\gamma)=\gamma\;,\;\gamma\in\mathfrak{Z}(E^{\sigma})$
follows easily.

Part (4) follows from (1) and (3). \end{proof}

Fix $X\in H^{\infty}(E)$ with $\Vert X\Vert\leq1$, let $\pi=\varphi_{\infty}\otimes I_{H}$,
as in Lemma \ref{1to1}, and let $\gamma$ be an element of $\mathbb{D}(\mathfrak{Z}(E^{\pi}))$.
Then if $\widehat{X}$ is the Schur class operator function on $\mathbb{D}((E^{\pi})^{\ast})$
determined by $X$ then by Corollary \ref{comp}, $\widehat{X}\circ g_{\gamma}$
also is a Schur class operator function on $\mathbb{D}((E^{\pi})^{\ast})$.
By Corollary~\ref{SchurMult} there is an element $\alpha_{\gamma}(X)$
in $H^{\infty}(E)$, whose norm does not exceed $1$, such that $\widehat{\alpha_{\gamma}(X)}=\widehat{X}\circ g_{\gamma}$.
Further, by Lemma~\ref{1to1}, this element is uniquely defined.
We can, of course, extend this to a map, $\alpha_{\gamma}$, from
$H^{\infty}(E)$ to itself such that, for $X\in H^{\infty}(E)$ and
$\eta\in\mathbb{D}((E^{\pi})^{\ast})$, \begin{equation}
\widehat{\alpha_{\gamma}(X)}(\eta^{\ast})=\widehat{X}(g_{\gamma}(\eta^{\ast})).\label{aut}\end{equation}

\begin{lemma} \label{autom} Let $\sigma$ and $\pi$ be as in Lemma
\ref{lifting}. Then:

\begin{enumerate}
\item [(i)] For every $\gamma\in\mathbb{D}(\mathfrak{Z}(E^{\pi}))$, $\alpha_{\gamma}$,
defined by equation (\ref{aut}) is an automorphism of the algebra
$H^{\infty}(E)$ that is completely isometric and is a homeomorphism
with respect to the ultraweak topology.
\item [(ii)] For every $\gamma\in\mathbb{D}(\mathfrak{Z}(E^{\sigma}))$
let $\alpha_{\gamma}$ be defined to be $\alpha_{\Phi(\gamma)}$ (with
$\Phi$ as in Lemma~\ref{lifting}). Then, 
for every $X\in H^{\infty}(E)$ and $\eta\in E^{\sigma}$, \begin{equation}
\widehat{\alpha_{\gamma}(X)}(\eta^{\ast})=\widehat{X}(g_{\gamma}(\eta^{\ast})).\label{autsig}\end{equation}

\end{enumerate}
\end{lemma}

\begin{proof} We first prove (i). Linearity and multiplicativity
of $\alpha_{\gamma}$ are easy to check. Since $g_{\gamma}^{2}=id$,
$\alpha_{\gamma}$ is invertible (with $\alpha_{\gamma}^{-1}=\alpha_{\gamma}$).
So it is an automorphism. Since $\alpha_{\gamma}$ maps the closed
unit ball of $H^{\infty}(E)$ into itself (as does the inverse map),
$\alpha_{\gamma}$ is isometric. It is, in fact, completely isometric.
To see this, consider, for $n\in\mathbb{N}$, the algebra $H^{\infty}(M_{n}(E))$,
associated with the $W^{*}$-correspondence $M_{n}(E)$ over the von
Neumann algebra $M_{n}(M)$. The corresponding Fock space is $M_{n}(\mathcal{F}(E))$
and the algebra can be identified with $M_{n}(H^{\infty}(E))$. The
representation $\sigma$ of $M$ gives rise to a representation $\sigma_{n}$
of $M_{n}(M)$ on $H^{(n)}=\mathbb{C}^{n}\otimes H$ (with $\sigma_{n}(M_{n}(M))^{\prime}=I_{\mathbb{C}^{n}}\otimes\sigma(M)^{\prime}\cong\sigma(M)^{\prime}$).
One can check that $E^{\sigma}\cong(M_{n}(E))^{\sigma_{n}}$. For
$\gamma\in\mathfrak{Z}(E^{\sigma})$, write $\gamma^{\prime}$ for
the corresponding element of $\mathfrak{Z}(M_{n}(E^{\sigma}))$. Then
$\alpha_{\gamma^{\prime}}$ acts on $M_{n}(H^{\infty}(E))$ by applying
$\alpha_{\gamma}$ to each entry. Since we know that $\alpha_{\gamma^{\prime}}$
is an isometry, it follows that $\alpha_{\gamma}$ is a complete isometry.

It is left to show that $\alpha_{\gamma}$ is continuous with respect
to the ultraweak topology.

For this, let $\{ X_{\beta}\}$ be a net in the closed unit ball of
$H^{\infty}(E)$ that converges ultraweakly to $X$. Since evaluating
at $\eta^{\ast}$ (for $\eta$ in the open unit ball) amounts to applying
a ultraweakly continuous representation , we have, for every such
$\eta$, $\widehat{X_{\beta}}(\eta^{\ast})\rightarrow\widehat{X}(\eta^{\ast})$
in the weak operator topology. Since this holds for $g_{\gamma}(\eta^{\ast})$
in place of $\eta$, we see that, for every $\eta$ in the open unit
ball of $E^{\sigma}$, \[
\widehat{\alpha_{\gamma}(X_{\beta})}(\eta^{\ast})\rightarrow\widehat{\alpha_{\gamma}(X)}(\eta^{\ast}).\]
 Using Lemma~\ref{1to1}, we find that $\alpha_{\gamma}(X_{\beta})\rightarrow\alpha_{\gamma}(X)$
in the ultraweak topology. This proves (i).

Part (ii) of the lemma results from the following computation \[
\widehat{\alpha_{\gamma}(X)}(\eta^{\ast})=\widehat{\alpha_{\Phi(\gamma)}(X)}(\Psi(\eta)^{\ast})|H=\widehat{X}(g_{\Phi(\gamma)}(\Psi(\eta)^{\ast}))|H\]
 \[
=\widehat{X}(g_{\Phi(\gamma)}(\Psi(\eta)^{\ast})|E\otimes H)=\widehat{X}(g_{\gamma}(\eta)^{\ast}),\]
 where we used equation (\ref{aut}) and Lemma~\ref{lifting}. \end{proof}

Note that we needed to use the representation $\pi$ in order to define,
for every $X\in H^{\infty}(E)$, the element $\alpha_{\gamma}(X)$
in $H^{\infty}(E)$ satisfying (\ref{aut}). That is, we used the
fact that $K(\pi)=0$. Once we defined it, it may be more convenient
to work with the original representation $\sigma$ (which can be chosen
to be an arbitrary faithful representation) and invoke (\ref{autsig}).
Note that, using Proposition~\ref{uniqueness}, we see that there
is only one automorphism that satisfies (\ref{autsig}).

\begin{theorem} \label{automcomp} Let $E$ be a $W^{\ast}$-correspondence
over $M$ and let $\sigma$ be a faithful normal representation of
$M$ on a Hilbert space $H$. Let $\alpha$ be an isometric automorphism
of $H^{\infty}(E)$ and assume that $g:\mathbb{D}(E^{\sigma})^{\ast}\rightarrow\mathbb{D}(E^{\sigma})^{\ast}$
is a biholomorphic automorphism of $\mathbb{D}(E^{\sigma})^{\ast}$
such that \[
\widehat{\alpha(X)}(\eta^{\ast})=\widehat{X}(g(\eta^{\ast}))\text{,}\]
 for all $X\in H^{\infty}(E)$ and all $\eta\in E^{\sigma}$. Then:

\begin{enumerate}
\item [(i)] $g(\mathbb{D}\mathfrak{Z}((E^{\sigma})^{\ast}))\subseteq\mathbb{D}\mathfrak{Z}((E^{\sigma})^{\ast})$. 
\end{enumerate}
\begin{itemize}
\item [(ii)] There is a $\gamma\in\mathbb{D}\mathfrak{Z}((E^{\sigma}))$
and a unitary operator $u$ in $\mathcal{L}(E)$ such that $u(\mathfrak{Z}(E))=\mathfrak{Z}(E)$
and such that \[
g(\eta^{\ast})=g_{\gamma}(\eta^{\ast})\circ(u\otimes I_{\mathcal{E}})\]
 (as a map from $E\otimes_{\sigma}H$ to $H$).
\item [(iii)] With $u$ as in (ii), there is an automorphism $\alpha_{u}$
of $H^{\infty}(E)$ such that $\alpha_{u}(T_{\xi})=T_{u\xi}$ for
every $\xi\in E$.
\item [(iv)] With $u$ and $\gamma$ as in (ii), \[
\alpha=\alpha_{\gamma}\circ\alpha_{u}\]
 where $\alpha_{\gamma}$ is the automorphism defined in equation
(\ref{aut}) (and satisfies (\ref{autsig})).
\item [(v)] For every $\eta_{1},\eta_{2},\ldots,\eta_{k}$ in the open
unit ball of $E^{\sigma}$, the map defined by the $k\times k$ matrix
\[
((id-\theta_{g(\eta_{i}^{\ast})^{\ast},g(\eta_{j}^{\ast})^{\ast}})\circ(id-\theta_{\eta_{i},\eta_{j}})^{-1})\]
 is completely positive. 
\end{itemize}
\end{theorem}

\begin{proof} Note first that, since $\alpha$ is an isometric automorphism,
it maps $\varphi_{\infty}(M)$ onto itself.

Suppose $\eta$ lies in $\mathbb{D}(\mathfrak{Z}(E^{\sigma})^{\ast})$.
Then, by part (3) of Lemma~\ref{propcenter}, $\widehat{X}(\eta^{\ast})\in\sigma(M)$
for every $X\in H^{\infty}(E)$. But then, for every $X$, $\widehat{X}(g(\eta^{\ast}))$
lies in $\sigma(M)$, showing that $g(\eta^{\ast})\in\mathfrak{Z}(E^{\sigma})$.
This proves (i).

The discussion following Lemma~\ref{propcenter} shows that we can
write $g=w\circ g_{\gamma}$ for some $\gamma$ in $\mathbb{D}\mathfrak{Z}((E^{\sigma}))$
and an isometry $w$ on $(E^{\sigma})^{\ast}$ that preserves the
center. Let $\alpha_{\gamma}$ be the automorphism described in Lemma~\ref{autom}(ii)
and write $\beta=\alpha_{\gamma}^{-1}\circ\alpha$. Then it follows
that \[
\widehat{\beta(X)}(\eta^{\ast})=\widehat{X}(w\eta^{\ast})\]
 for $X\in H^{\infty}(E)$ and $\eta\in\mathbb{D}(E^{\sigma})$.

For $\eta=0$ and $Y\in H^{\infty}(E)$ we have $\widehat{Y}(0)=\sigma(\mathbb{E}_{0}(Y))$
where $\mathbb{E}_{0}$ is the conditional expectation of $H^{\infty}(E)$
onto $M$ (where $M$ is viewed as the {}``zeroth term\textquotedblright).
Thus, $\sigma(\mathbb{E}_{0}(\beta(X)))=\widehat{\beta(X)}(0)=\widehat{X}(0)=\sigma(\mathbb{E}_{0}(X))$
for every $X\in H^{\infty}(E)$. Since $\sigma$ is faithful, $\mathbb{E}_{0}(\beta(X))=\mathbb{E}_{0}(X)$.
Thus, for every $\xi\in E$, $\mathbb{E}_{0}(\beta(T_{\xi}))=0$ and
we can write \begin{equation}
\beta(T_{\xi})=T_{\theta}+Y\label{beta}\end{equation}
 where $Y$ lies in $(T_{E})^{2}H^{\infty}(E)$. Write $C$ for $(T_{E})^{2}H^{\infty}(E)$.
Since (\ref{beta}) holds for all $\xi\in E$, $\beta(C)\subseteq C$.
We can apply the same arguments to $\beta^{-1}$, in place of $\beta$,
and find that $\beta^{-1}(C)\subseteq C$. Applying $\beta^{-1}$
to (\ref{beta}), we find that \begin{equation}
\beta^{-1}(T_{\theta})=T_{\xi}+Z\label{betainv}\end{equation}
 for some $Z\in C$.

Arguing as in the proof of Proposition~\ref{uniqueness}, we find
that, for every $h\in H$, $\Vert\beta(T_{\xi})h\Vert\leq\Vert T_{\xi}h\Vert$
and $\Vert\beta(T_{\xi})\Vert^{2}=\Vert Yh\Vert^{2}+\Vert T_{\theta}h\Vert^{2}\geq\Vert T_{\theta}h\Vert^{2}$.
Thus $\Vert T_{\xi}h\Vert\geq\Vert T_{\theta}h\Vert$. Applying the
same arguments to $\beta^{-1}$ (using (\ref{betainv}) in place of
(\ref{beta})) we find that $\Vert T_{\theta}h\Vert\geq\Vert T_{\xi}h\Vert$
and, thus, $\Vert T_{\xi}h\Vert=\Vert T_{\theta}h\Vert$ and, consequently,
$Yh=0$ for all $h\in H$. Thus $Y=0$ and 
$\beta(T_{\xi})=T_{\theta}$. Since $\beta$ is isometric, $\Vert T_{\xi}\Vert=\Vert T_{\theta}\Vert$.
It follows that $\Vert\xi\Vert=\Vert\theta\Vert$. If we write $\theta=u\xi$
(and recall that then $\beta(T_{\xi})=T_{u\xi}$) then $u$ is a linear
isometry. We also have, for $a\in M$, $T_{u(\xi a)}=\beta(T_{\xi a})=\beta(T_{\xi}a)=\beta(T_{\xi})a=T_{u(\xi)}a=T_{u(\xi)a}$.
Hence $u$ is an isometric (right) module map and, therefore, $u$
lies in $\mathcal{L}(E)$. Since $\beta$ is an automorphism, $u$
is a unitary operator. We also have $\beta(T_{\xi})=T_{u\xi}$, so
$\beta=\alpha_{u}$ (in the notation of (iii)). This proves (iii)
and (iv).

Recall that $\widehat{\beta(X)}(\eta^{\ast})=\widehat{X}(w\eta^{\ast})$
and set $X=T_{\xi}$ to get $\widehat{T_{u\xi}}(\eta^{\ast})=\widehat{\beta(T_{\xi})}(\eta^{\ast})=\widehat{T_{\xi}}(w\eta^{\ast})$.
Hence $\eta^{\ast}L_{u\xi}=(w\eta^{\ast})L_{\xi}$. Applying this
to $h\in\mathcal{E}$ we get $\eta^{\ast}(u\xi\otimes h)=(w\eta^{\ast})(\xi\otimes h)$.
Hence $w\eta^{\ast}=\eta^{\ast}\circ(u\otimes I)$, proving $g(\eta^{\ast})=g_{\gamma}(\eta^{\ast})\circ(u\otimes I_{\mathcal{E}})$.
To prove (ii) we need only to show that $u$ preserves the center
of $E$. So fix $\xi\in\mathfrak{Z}(E)$. By Lemma~\ref{propcenter},
$L_{\xi}^{\ast}$ lies in the center of $(E^{\sigma})^{\ast}$. Thus
$wL_{\xi}^{\ast}$ lies in $\mathfrak{Z}((E^{\sigma})^{\ast})$. But
$wL_{\xi}^{\ast}=L_{\xi}^{\ast}\circ(u\otimes I)=L_{u^{\ast}\xi}$.
Thus $L_{u^{\ast}\xi}$ lies in $\mathfrak{Z}((E^{\tau})^{\ast})$.
Using Lemma~\ref{propcenter} again we get $u^{\ast}\xi\in\mathfrak{Z}(E)$.
This shows that $u^{\ast}\mathfrak{Z}(E)\subseteq\mathfrak{Z}(E)$
and, applying the same argument to $\beta^{-1}$, we complete the
proof of (ii).

To prove (v), fix $b\in\sigma(M)^{\prime}$ and $\eta_{i},\eta_{j}$
in $\mathbb{D}(E^{\sigma})$ and compute $\langle g(\eta_{i}^{\ast}),b\cdot g(\eta_{j}^{\ast})\rangle=g(\eta_{i}^{\ast})(I_{E}\otimes b)g(\eta_{j}^{\ast})^{\ast}=g_{\gamma}(\eta_{i}^{\ast})(u\otimes I_{\mathcal{E}})(I_{E}\otimes b)(u^{\ast}\otimes I_{\mathcal{E}})g_{\gamma}(\eta_{j})^{\ast}=g_{\gamma}(\eta_{i}^{\ast})(I_{E}\otimes b)g_{\gamma}(\eta_{j}^{\ast})^{\ast}=\langle g_{\gamma}(\eta_{i}^{\ast}),b\cdot g_{\gamma}(\eta_{j}^{\ast})\rangle$.
Thus (v) follows from Corollary~\ref{kerg}. \end{proof}

Combining Theorem~\ref{automcomp} with Theorem~\ref{autom1}, we
get the following.

\begin{theorem} \label{comb}Let $E$ be a faithful $W^{\ast}$-correspondence
over $M$ where $\mathfrak{Z}(M)$ is atomic. Let $\alpha$ be an
automorphism of $H^{\infty}(E)$ that is completely isometric and
a $w^{\ast}$-homeomorphism and leaves $\varphi_{\infty}(M)$ elementwise
fixed and let $\sigma$ be a faithful representation of $M$.

Then there is a $\gamma\in\mathbb{D}\mathfrak{Z}((E^{\sigma}))$ and
a unitary operator $u$ in $\mathcal{L}(E)$, satisfying $u(\mathfrak{Z}(E))=\mathfrak{Z}(E)$,
such that \[
\alpha=\alpha_{\gamma}\circ\alpha_{u},\]
 where $\alpha_{\gamma}$ is the automorphism defined in Lemma~\ref{autom}
and $\alpha_{u}(T_{\xi})=T_{u\xi}$ for every $\xi\in E$.

In particular, if $\mathfrak{Z}(E)=\{0\}$, every such automorphism
is $\alpha_{u}$ for some unitary operator $u\in\mathcal{L}(E)$.
\end{theorem}

Theorem \ref{comb} provides another perspective on the results from
\cite{MS86,MS87}. The analytic crossed products discussed there are
of the form $H^{\infty}(E)$, where $E$ is the correspondence $_{\alpha}M$
associated with a von Neumann algebra $M$ and an automorphism $\alpha$
that is properly outer. This means that $\mathfrak{Z}(E)=\{0\}$.
Theorem \ref{comb} implies that all automorphisms of $H^{\infty}(E)$
are given by automorphisms of $\dot{M}$.

\section{Examples : Graph Algebras}

In this section we consider some examples that come from directed
graphs. We shall assume for simplicity that our graphs have finitely
many vertices and edges. We write $\mathcal{Q}$ both for the graph
and for its set of edges. \ The space of vertices will be denoted
$V$. We shall write $s$ and $r$ for the source and range maps on
$\mathcal{Q}$, mapping $\mathcal{Q}$ to $V$, and we shall think
of an edge $e$ in $\mathcal{Q}$ as {}``pointing\textquotedblright\ from
$s(e)$ to $r(e)$. For simplicity, we shall also assume that $r$
is surjective, i.e., we shall assume that $\mathcal{Q}$ is \emph{without
sources}. Write $\mathcal{Q}^{\ast}$ for the set of all finite paths
in $\mathcal{Q}$, i.e., the path category generated by $\mathcal{Q}$.
An element in $\mathcal{Q}$ will be written $\alpha=e_{1}e_{2}\cdots e_{k}$,
where $s(e_{i})=r(e_{i+1})$. We set $s(\alpha)=s(e_{k})$, $r(\alpha)=r(e_{1})$,
and $|\alpha|=k,$ the length of $\alpha$. We will also view vertex
$v\in V$ as a {}``path of length $0$\char`\"{}, and we extend $r$
and $s$ to $V$ simply by setting $r(v)=s(v)=v.$

Let $M$ be $C(V),$ the set of complex-valued functions on $V$.
Of course, $M$ is a finite dimensional commutative von Neumann algebra.
Likewise, we let $E$ be $C(\mathcal{Q})$, the set of complex-valued
functions on $\mathcal{Q}$. Then we define an $M$-bimodule structure
on $E$ as follows: for $f\in E$, $\psi\in M$ and $e\in\mathcal{Q}$,
\[
(f\psi)(e):=f(e)\psi(s(e)),\]
 and \[
(\psi f)(e):=\psi(r(e))f(e).\]
 Note that the {}``no sources\char`\"{} assumption implies that the
left action of $M$ is faithful. An $M$-valued inner product on $E$
will be given by the formula \[
\langle f,g\rangle(v)=\sum_{s(e)=v}\overline{f(e)}g(e),\]
 for $f,g\in E$ and $v\in V$. With these operations, $E$ becomes
a $W^{\ast}$-correspondence over $M$. The algebra $H^{\infty}(E)$
in this case will be written $H^{\infty}(\mathcal{Q})$. In the literature,
$H^{\infty}(\mathcal{Q})$ is sometimes denoted $\mathcal{L}_{\mathcal{Q}}$.
It is the ultraweak closure of the tensor algebra $\mathcal{T}_{+}(E(\mathcal{Q}))$
acting on the Fock space of $\mathcal{F}(E(\mathcal{Q}))$. For $e\in\mathcal{Q}$,
let $\delta_{e}$ be the $\delta$-function at $e$, i.e., $\delta_{e}(e^{\prime})=1$
if $e=e^{\prime}$ and is zero otherwise. Then $T_{\delta_{e}}$ is
a partial isometry that we denote by $S_{e}$. Also, for $v\in V$,
$P_{v}$ is defined to be $\varphi_{\infty}(\delta_{v})$. Then each
$P_{v}$ is a projection and it is an easy matter to see that the
families $\{ S_{e}:e\in\mathcal{Q}\}$ and $\{ P_{v}:v\in V\}$ form
a \emph{Cuntz-Toeplitz family} in the sense that the following conditions
are satisfied:

\begin{enumerate}
\item [(i)] $P_{v}P_{u}=0$ if $u\neq v$,
\item [(ii)] $S_{e}^{*}S_{f}=0$ if $e\neq f$
\item [(iii)] $S_{e}^{*}S_{e}=P_{s(e)}$ and
\item [(iv)] $\sum_{r(e)=v}S_{e}S_{e}^{*}\leq P_{v}$ for all $v\in V$. 
\end{enumerate}
In fact, these particular families yield a faithful representation
of the Cuntz-Toeplitz algebra $\mathcal{T}(E(\mathcal{Q}))$ \cite{FR99}.
The algebra $\mathcal{T}_{+}(E(\mathcal{Q}))$ is the norm-closed
(unstarred) algebra that they generate inside $\mathcal{T}(E(\mathcal{Q}))$
and $H^{\infty}(\mathcal{Q})$ is the ultraweak closure of $\mathcal{T}_{+}(E(\mathcal{Q}))$.
The algebra $\mathcal{T}_{+}(E(\mathcal{Q}))$ was first defined and
studied in \cite{pM97}, providing examples of the theory developed
in \cite{MS98}. It was called a quiver algebra there because in pure
algebra, graphs of the form $\mathcal{Q}$ are called quivers. (Hence
the notation we use here.) The properties of quiver algebras were
further developed in \cite{MS99}. In \cite{KP}, the focus was on
$H^{\infty}(\mathcal{Q})$ and the authors called this algebra a free
semigroupoid algebras. Both algebras are often represented as algebras
of operators on $l_{2}(\mathcal{Q}^{\ast})$, and it will be helpful
to understand how from the perspective of this note. Let $H_{0}$
be a Hilbert space whose dimension equals the number of vertices,
let $\{ e_{v}|\; v\in V\}$ be a fixed orthonormal basis for $H_{0}$
and let $\sigma_{0}$ be the diagonal representation of $M=C(V)$
on $H_{0}$. Then $l_{2}(\mathcal{Q}^{\ast})$ is isomorphic to $\mathcal{F}(E(\mathcal{Q}))\otimes_{\sigma_{0}}H_{0}$
where the isomorphism maps an element $\xi_{\alpha}$ of the standard
orthonormal basis of $l_{2}(\mathcal{Q}^{\ast})$ to $\delta_{\alpha}\otimes e_{s(e)}$
(where, for $\alpha=e_{1}\cdots e_{k}$, $\delta_{\alpha}=\delta_{e_{1}}\otimes\cdots\otimes\delta_{e_{k}}\in E^{\otimes k}$).
The partial isometries $S_{e}$ can then be viewed as the shift operators
$S_{e}\xi_{\alpha}=\xi_{e\alpha}$. Thus, the representations of $\mathcal{T}_{+}(E(\mathcal{Q}))$
and $H^{\infty}(\mathcal{Q})$ on $l_{2}(\mathcal{Q}^{\ast})$ are
just the representations induced by $\sigma_{0}$.

Quite generally, a completely contractive covariant representation
of $E(\mathcal{Q})$ on a Hilbert space $H$ is given by a representation
$\sigma$ of $M=C(V)$ on $H$ and by a contractive map $\tilde{T}:E\otimes_{\sigma}H\rightarrow H$
satisfying equation (\ref{covariance}). The representation $\sigma$
is given by the projections $Q_{v}=\sigma(\delta_{v})$ whose sum
is $I$. Also, from $\tilde{T}$ we may define maps $T(e)\in B(H)$
by the equation $T(e)h=\tilde{T}(\delta_{e}\otimes h)$ and it is
easy to check that $\tilde{T}\tilde{T}^{\ast}=\sum_{e}T(e)T(e)^{\ast}$
and $T(e)=Q_{r(e)}T(e)Q_{s(e)}$. Thus to every completely contractive
representation of the quiver algebra $\mathcal{T}_{+}(E(\mathcal{Q}))$
we associate a family $\{ T(e)|e\in\mathcal{Q}\}$ of maps on $H$
that satisfy $\sum_{e}T(e)T(e)^{\ast}\leq I$ and $T(e)=Q_{r(e)}T(e)Q_{s(e)}$.
Conversely, every such family defines a representation, written $\sigma\times T$
(or $\sigma\times\tilde{T}$), satisfying $(\sigma\times T)(S_{e})=T(e)$
and $(\sigma\times T)(P_{v})=Q_{v}$.

We fix $\sigma$ to be $\sigma_{0}$ and write $H$ in place of $H_{0}$.
So that, in this case, each projection $Q_{v}$ is one dimensional
(with range equal to $\mathbb{C}e_{v}$). Then obviously $\sigma(M)^{\prime}=\sigma(M)$.
To describe the $\sigma$-dual of $E$, write $\mathcal{Q}^{-1}$
for the directed graph obtained from $\mathcal{Q}$ by reversing all
arrows, so that $s(e^{-1})=r(e)$ and $r(e^{-1})=s(e)$. Sometimes
$\mathcal{Q}^{-1}$ is denoted $\mathcal{Q}^{op}$ and is called the
opposite graph. Note that the Hilbert space $E\otimes_{\sigma}H_{0}$
is spanned by the orthonormal basis $\{\delta_{e}\otimes e_{s(\alpha)}\}$.
Fix $\eta\in E^{\sigma}$ and note that its covariance property implies
that, for every $e\in\mathcal{Q}$, $\eta^{\ast}(\delta_{e}\otimes e_{s(e)})=\eta^{\ast}(\delta_{r(e)}\delta_{e}\otimes e_{s(e)})=Q_{r(e)}\eta^{\ast}(\delta_{e}\otimes e_{s(e)})=\overline{\eta(e^{-1})}e_{r(e)}$
for some $\overline{\eta(e^{-1})}\in\mathbb{C}$. The reason for the
{}``strange\char`\"{} way of writing that scalar is that we can view
$\eta$ as an element of $E(\mathcal{Q}^{-1})$ and the correspondence
structure on $E^{\sigma},$ as described in Proposition~\ref{corres},
fits the correspondence structure of $E(\mathcal{Q}^{-1})$. Consequently,
we can identify the two and write \[
E^{\sigma}=E(\mathcal{Q}^{-1}).\]
 (See Example 4.3 in \cite{MS03} for a description of the structure
of the dual correspondence for more general representations $\sigma$
). It will also be convenient to write $\eta$ matricially with respect
to the orthonormal bases $\{\delta_{v}\mid v\in V\}$ of $H_{0}$
and $\{\delta_{e}\otimes e_{s(e)}\}_{e\in\mathcal{Q}}$ of $E\otimes H_{0}$
as \begin{equation}
(\eta)_{e,r(e)}=\eta(e^{-1}).\label{matrix1}\end{equation}

Suppose $\eta\in\mathbb{D}(E^{\sigma})$. For every $X\in H^{\infty}(\mathcal{Q})$,
we have defined $X(\eta^{\ast})$ as an element of $B(H)$ in Remark~\ref{evaluation}.
For the generators of $H^{\infty}(\mathcal{Q})$, the definition yields
the equations, \begin{equation}
\widehat{P_{v}}(\eta^{\ast})=\theta_{v,v}\;,\; v\in V\label{eval1}\end{equation}
 and \begin{equation}
\widehat{S_{e}}(\eta^{\ast})=\overline{\eta(e^{-1})}\theta_{r(e),s(e)}\;,\; e\in\mathcal{Q}\label{eval2}\end{equation}
 where $\theta_{v,w}$ is the partial isometry operator on $H$ that
maps $e_{w}$ to $e_{v}$ and vanishes on $(e_{w})^{\perp}$. For
a general $X\in H^{\infty}(\mathcal{Q})$, $\widehat{X}(\eta^{\ast})$
is obtained by using the linearity, multiplicativity and $w^{\ast}$-continuity
of the map $X\mapsto\widehat{X}(\eta^{\ast})$.

The proof of the next lemma is straightforward and is omitted.

\begin{lemma} \label{centerq} The centers of the correspondences
$E(\mathcal{Q})$ and $E(\mathcal{Q}^{-1})$ are given by the formulae
\[
\mathcal{\mathfrak{Z}}(E(\mathcal{Q}))=span\{\delta_{e}\;|\; s(e)=r(e)\}\]
 and \[
\mathcal{\mathfrak{Z}}(E(\mathcal{Q}^{-1}))=span\{\delta_{e^{-1}}\;|\; s(e)=r(e)\}.\]

\end{lemma}

The following proposition is immediate from Theorem~\ref{comb}.

\begin{proposition} \label{noloops} If there is no $e\in\mathcal{Q}$
with $s(e)=r(e)$, then every automorphism $\alpha$ of $H^{\infty}(\mathcal{Q})$
that is completely isometric, $w^{\ast}$-homeomorphic and leaves
$\varphi_{\infty}(C(V))$ elementwise fixed (that is, does not permute
the vertices) is of the form $\alpha_{u}$ for some unitary $u\in\mathcal{L}(E(\mathcal{Q}))$.
That is, \[
\alpha(S_{e})=\sum_{s(f)=s(e)}u_{f,e}S_{f}\]
 where the scalars $u_{f,e}$ are given by $u_{f,e}=(u(\delta_{e}))(f)$.
(Note that this is zero if $s(f)\neq s(e)$, since $u(\delta_{e})=u(\delta_{e}\delta_{s(e)})=u(\delta_{e})\delta_{s(e)}$).
\end{proposition}

We note, as we did at the beginning of Section 4, that the assumptions
made on the automorphism can be weakened using arguments of \cite{KK04}
but we shall not elaborate on this here.

\begin{example} \label{cycle} Let $\mathcal{Q}$ be an n-cycle (for
$n>1$) ; that is $V=\{ v_{1},v_{2},\ldots,v_{n}\}$ and $\mathcal{Q}=\{ e_{1},\ldots,e_{n}\}$
where $e_{i}$ is the arrow from $v_{1}$ to $v_{i+1}$ (or to $v_{1}$
when $i=n$). Then, for every $\alpha$ as in Proposition~\ref{noloops},
there are $\{\lambda_{1},\lambda_{2},\ldots,\lambda_{n}\}$ with $|\lambda_{i}|=1$,
such that $\alpha(S_{e_{i}})=\lambda_{i}S_{e_{i}}$ for all $i$.
\end{example}

The rest of this section will be devoted to the study of the following
example, which is very simple, yet provides a full array of the structures
we have been studying.

\begin{example} Let the vertex set of the graph have two elements:
$V=\{ v,w\}.$ Suppose the edge set consists of three elements $\mathcal{Q}=\{ e,f,g\}$,
where $e$ is the arrow from $v$ to $w$, so $s(e)=v$, $r(e)=w$;
$f$ is an arrow from $w$ to $v;$ and $g$ is a loop based at $w$,
$s(g)=r(g)=w$. \end{example}

Then by Lemma \ref{centerq}, $\mathcal{\mathfrak{Z}}(E(\mathcal{Q)})=\mathbb{C}\delta_{g}$.
We know from Theorem~\ref{comb} that every automorphism $\alpha$
is the composition of an automorphism, written $\alpha_{u}$ associated
with a unitary in $\mathcal{L}(E(\mathcal{Q}))$ that maps $\delta_{g}$
into $\lambda_{3}\delta_{g}$ (with $|\lambda_{3}|=1$) and an automorphism
associated with a {}``Möbius transformation\char`\"{}.

As noted in Proposition~\ref{noloops}, $(u(\delta_{e^{\prime}}))(f^{\prime})=0$
unless $s(e^{\prime})=s(f^{\prime})$, so that $u(\delta_{e})\in\mathbb{C}\delta_{e}$
and $u(\delta_{f})\in span\{\delta_{f},\delta_{g}\}$. Since $u^{*}$
is unitary, we have that $u(\delta_{f})=\lambda_{f}\delta_{f}$. Thus
\begin{equation}
\alpha_{u}(S_{e})=\lambda_{e}S_{e},\;\;\alpha_{u}(S_{f})=\lambda_{f}S_{f}\label{u}\end{equation}
 and \[
\alpha_{u}(S_{g})=\lambda_{g}S_{g}\]
 for $\lambda_{e},\lambda_{f},\lambda_{g}$ with absolute value $1$.

It is left to analyze the Möbius transformations and the corresponding
automorphisms. Since the center of $E^{\sigma}$ are scalar multiples
of $\delta_{g^{-1}}$, the Möbius transformations are associated with
scalars $\lambda\in\mathbb{D}$ (in fact, with $\lambda\delta_{g^{-1}}$)
and will be denoted $\tau_{\lambda}$, $\lambda\in\mathbb{D}$. We
have \begin{equation}
\tau_{\lambda}(\eta^{\ast})=\Delta_{\lambda}(I-\eta^{\ast}(\lambda\delta_{g^{-1}}))^{-1}(\bar{\lambda}\delta_{g^{-1}}-\eta^{\ast})\Delta_{\lambda\ast}^{-1}\label{taulambda}\end{equation}
 where $\Delta_{\lambda}=(I_{H}-(\lambda\delta_{g^{-1}})^{\ast}(\lambda\delta_{g^{-1}}))^{1/2}$
and $\Delta_{\lambda\ast}=(I_{E\otimes H}-(\lambda\delta_{g^{-1}})(\lambda\delta_{g^{-1}})^{\ast})^{1/2}$.
It will be convenient to write $\tau_{\lambda}(\eta^{\ast})$ matricially
as a map from $E\otimes H$, with the ordered orthonormal basis $\{\delta_{e}\otimes\delta_{v},\delta_{f}\otimes\delta_{w},\delta_{g}\otimes\delta_{w}\}$,
to $H$, with the ordered orthonormal basis $\{\delta_{v},\delta_{w}\}$.
Using the formula (\ref{matrix1}), we see that \[
\eta=\left(\begin{array}{cc}
0 & \eta(e^{-1})\\
\eta(f^{-1}) & 0\\
0 & \eta(g^{-1})\end{array}\right)\]
 and \[
\lambda\delta_{g^{-1}}=\left(\begin{array}{cc}
0 & 0\\
0 & 0\\
0 & \lambda\end{array}\right).\]

The computation of the expression in (\ref{taulambda}) yields \[
\tau_{\lambda}(\eta^{*})=\left(\begin{array}{ccc}
0 & -\overline{\eta(f^{-1})} & 0\\
\frac{-\overline{\eta(e^{-1})}(1-|\lambda|^{2})^{1/2}}{1-\lambda\overline{\eta(g^{-1})}} & 0 & \frac{\bar{\lambda}-\overline{\eta(g^{-1})}}{1-\lambda\overline{\eta(g^{-1})}}\end{array}\right).\]
 Thus \[
\overline{\tau_{\lambda}(\eta^{*})^{*}(e^{-1})}=\frac{-\overline{\eta(e^{-1})}(1-|\lambda|^{2})^{1/2}}{1-\lambda\overline{\eta(g^{-1})}}=-\overline{\eta(e^{-1})}(1-|\lambda|^{2})^{1/2}\sum_{k=0}^{\infty}(\lambda\overline{\eta(g^{-1})})^{k},\]
 \[
\overline{\tau_{\lambda}(\eta^{*})^{*}(f^{-1})}=-\overline{\eta(f^{-1})},\]
 and \[
\overline{\tau_{\lambda}(\eta^{*})^{*}(g^{-1})}=\frac{\bar{\lambda}-\overline{\eta(g^{-1})}}{1-\lambda\overline{\eta(g^{-1})}}=(\bar{\lambda}-\overline{\eta(g^{-1})})\sum_{k=0}^{\infty}(\lambda\overline{\eta(g^{-1})})^{k}.\]

This suggests setting \[
T(e)=-(1-|\lambda|^{2})^{1/2}\sum_{k=0}^{\infty}(\lambda S_{g})^{k}S_{e},\]
 \[
T(f)=-S_{f}\]
 and \[
T(g)=-(\bar{\lambda}P_{w}-S_{g})\sum_{k=0}^{\infty}(\lambda S_{g})^{k}.\]

Using (\ref{eval1}), (\ref{eval2}) and the fact that the map $X\mapsto\widehat{X}(\eta^{*})$
is a continuous homomorphism, we get \[
\widehat{T(e)}(\eta^{*})=\overline{\tau_{\lambda}(\eta^{*})^{*}(e^{-1})}\theta_{w,v}\]
 , \[
\widehat{T(f)}(\eta^{*})=\overline{\tau_{\lambda}(\eta^{*})^{*}(f^{-1})}\theta_{v,w}\]
 and \[
\widehat{T(g)}(\eta^{*})=\overline{\tau_{\lambda}(\eta^{*})^{*}(g^{-1})}\theta_{w,w}.\]

Using Theorem~\ref{autom1}, Theorem~\ref{comb}, Equation (\ref{u})
and Theorem~\ref{uniqueness}, we conclude the following.

\begin{theorem} \label{automgraph}

\begin{enumerate}
\item [(1)] For every $\lambda\in\mathbb{D}$, there is a unique automorphism
$\alpha_{\lambda}$ of $H^{\infty}(\mathcal{Q})$ such that, for every
$e^{\prime}\in\{ e,f,g\}$, $\alpha_{\lambda}(S_{e^{\prime}})-T(e^{\prime})\in K(\sigma)$.
\item [(2)] Every completely isometric, $w^{\ast}$-homeomorphic auto\-mor\-phism
$\alpha$ of $H^{\infty}(\mathcal{Q})$ can be written \[
\alpha=\alpha_{u}\circ\alpha_{\lambda}\]
 where $\lambda\in\mathbb{D}$ and $\alpha_{u}(S_{e^{\prime}})=\lambda_{e^{\prime}}S_{e^{\prime}}$
for every $e^{\prime}\in\{ e,f,g\}$ (where $\lambda_{e},\lambda_{f}$
and $\lambda_{g}$ are complex numbers of absolute value $1$). 
\end{enumerate}
\end{theorem}

\begin{proof} The only thing that we need to clarify here is that,
in part (2), we do not have to require that $\alpha$ fixes $P_{v}$
and $P_{w}$. Indeed, assume that $\alpha$ satisfies $\alpha(P_{v})=P_{w}$
and $\alpha(P_{w})=P_{v}$. Then $\alpha(S_{e})=P_{v}\alpha(S_{e})P_{w}$
and, thus, $\mathbb{E}_{0}(\alpha(S_{e}))=0$ and $\mathbb{E}_{1}(\alpha(S_{e}))\in\mathbb{C}S_{f}$.
Similarly, we get $\mathbb{E}_{0}(\alpha(S_{f}))=\mathbb{E}_{1}(\alpha(S_{g}))=0$,
$\mathbb{E}_{1}(\alpha(S_{f}))\in\mathbb{C}S_{e}$ and $\mathbb{E}_{0}(\alpha(S_{g}))\in\mathbb{C}P_{v}$.
Thus, $S_{g}$ is not in the range of $\alpha$, contradicting the
surjectivity of $\alpha$. \end{proof}

Finally, we note the following.

\begin{proposition} \label{K} In this example, $K(\sigma)$ is the
ideal generated by the commutator $[S_{g},S_{e}S_{f}]$. \end{proposition}

\begin{proof} Since we shall not use this result, we only sketch
the idea of the proof. It follows from Lemma~\ref{kersigma} that
it suffices to analyze $\mathbb{E}_{k}(K(\sigma))$ for a given $k$.
Since $K(\sigma)$ is an ideal, it suffices to consider $P_{v^{\prime}}\mathbb{E}_{k}(K(\sigma))P_{v^{\prime\prime}}$
for fixed $v^{\prime},v^{\prime\prime}\in\{ v,w\}$. Evaluating an
element of $P_{v^{\prime}}\mathbb{E}_{k}(K(\sigma))P_{v^{\prime\prime}}$
in $\eta^{\ast}$ yields a polynomial in three the variables $z_{1}=\overline{\eta(e^{-1})},z_{2}=\overline{\eta(f^{-1})}$
and $z_{3}=\overline{\eta(f^{-1})}$. This polynomial is defined on
a small enough neighborhood of $0$ and, from the definition of $K(\sigma)$,
it vanishes there. It follows that its coefficients are all $0$.
This shows that an element in $P_{v^{\prime}}\mathbb{E}_{k}(K(\sigma))P_{v^{\prime\prime}}$
is a linear combination of sums of the form $\sum a_{i}S_{\alpha_{i}}$
(for some paths $\alpha_{i}$) where $\sum a_{i}=0$ and for every
$i,j$, the paths $\alpha_{i}$ and $\alpha_{j}$ satisfy $s(\alpha_{i})=s(\alpha_{j})=v^{\prime\prime}$,
$r(\alpha_{i})=r(\alpha_{j})=v^{\prime}$ and both paths contain the
same edges (with the same multiplicities) but in a different order.
A moment's reflection shows that this can happen only if the two paths
are identical except that, at some points, one path travels along
$g$ and then along $ef$ while the other path {}``chooses\char`\"{}
to travel first along $ef$ and then along $g$. This shows that the
element in $P_{v^{\prime}}\mathbb{E}_{k}(K(\sigma))P_{v^{\prime\prime}}$
lies in the ideal generated by $[S_{g},S_{e}S_{f}]$. \end{proof}

\end{document}